\documentclass[11pt]{article}

\usepackage{bbm}

\usepackage{color}
\usepackage{amsmath}
\usepackage{latexsym,amssymb}
\usepackage{graphicx}
\usepackage{stmaryrd}
\usepackage{fixmath}

\newtheorem{theorem}{Theorem}[section]
\newtheorem{lemma}[theorem]{Lemma}

\newtheorem{proposition}[theorem]{Proposition}

\newtheorem{remark}[theorem]{Remark}

\numberwithin{equation}{section}

\newcommand{\noi}{\noindent}

\newcommand{\E}{\mathbb{E}}
\newcommand{\R}{\mathbb{R}}
\newcommand{\N}{\mathbb{N}}

\newcommand{\la}{\lambda}
\newcommand{\sig}{\sigma}

\newcommand{\eps}{\varepsilon}
\newcommand{\ph}{\varphi}

\newcommand{\kap}{\kappa}

\newcommand{\del}{\delta}

\newcommand{\Gam}{\mathnormal{\Gamma}}

\newcommand{\PI}{\mathnormal{\Pi}}

\newcommand{\Om}{\mathnormal{\Omega}}

\newcommand{\Z}{{\mathbb Z}}

\newcommand{\PP}{{\mathbb P}}

\newcommand{\bE}{{\mathbf E}}
\newcommand{\bP}{{\mathbf P}}
\newcommand{\bQ}{{\mathbf Q}}

\newcommand{\calA}{{\cal A}}

\newcommand{\calC}{{\cal C}}
\newcommand{\calD}{{\cal D}}

\newcommand{\calF}{{\cal F}}
\newcommand{\calG}{{\cal G}}

\newcommand{\calK}{{\cal K}}
\newcommand{\calL}{{\cal L}}
\newcommand{\calM}{{\cal M}}

\newcommand{\calS}{{\mathcal{S}}}

\newcommand{\calU}{{\cal U}}

\newcommand{\skp}{\vspace{\baselineskip}}

\newcommand{\diag}{{\rm diag}}

\newcommand{\w}{\wedge}

\newcommand{\To}{\Rightarrow}

\newcommand{\dist}{{\rm dist}}
\newcommand{\iy}{\infty}
\newcommand{\ds}{\displaystyle}

\newcommand{\down}{\downarrow}

\newcommand{\qed}{\hfill $\Box$}

\newcommand{\osc}{{\rm osc}}
\newcommand{\one}{\mathbbm{1}}

\newcommand{\Erbm}{\E}
\newcommand{\PPrbm}{\PP}

\begin{document}

\title{Serve the shortest queue and Walsh Brownian motion\thanks{This is the final version of the paper. To appear in {\it The Annals of Applied Probability.}}}

\author{Rami Atar\thanks{Department of Electrical Engineering,
Technion--Israel Institute of Technology, Haifa 32000, Israel. Research supported in part by the ISF (grant 1184/16)}
\and
Asaf Cohen\thanks{Department of Statistics,
	University of Haifa,
	Haifa, 31905, Israel,
	Email:
	shloshim@gmail.com,
	web: https://sites.google.com/site/asafcohentau/
}}

\date{July 15, 2018}

\maketitle

\begin{abstract}
We study a single-server Markovian queueing model with $N$ customer classes in which
priority is given to the shortest queue.
Under a critical load condition,
we establish the diffusion limit of the nominal workload and queue length processes
in the form of a Walsh Brownian motion (WBM)
living in the union of the $N$ nonnegative coordinate axes in $\R^N$
and a linear transformation thereof.
This reveals the following asymptotic behavior.
Each time that queues begin to build starting from an empty system,
one of them becomes dominant in the sense that
it contains nearly all the workload in the system,
and it remains so until the system becomes (nearly) empty again.
The radial part of the WBM, given as a reflected Brownian motion (RBM)
on the half-line, captures the total workload asymptotics,
whereas its angular distribution expresses how likely it is for each
class to become dominant on excursions.

As a heavy traffic result it is nonstandard in three ways:
(i) In the terminology of Harrison \cite{har-bal} it is {\it unconventional},
in that the limit is not an RBM.
(ii) It does not constitute an {\it invariance principle},
in that the limit law (specifically, the angular distribution)
is not determined solely by the first two moments of the data,
and is sensitive even to tie breaking rules.
(iii) The proof method does not fully characterize
the limit law (specifically, it gives no information on the angular distribution).

\skp

\noi{\bf AMS subject classification:}
60F05, 
93E03, 
60K25, 
60J65, 
60J70

\skp

\noi{\bf Keywords:}
Serve the shortest queue, heavy traffic, diffusion limits,
Walsh Brownian motion

\end{abstract}

\section{Introduction}\label{sec1}

We consider a multiclass single-server queueing system operating under
{\it serve the shortest queue} (SSQ) (also referred to in the literature as
{\it shortest queue first}) regime,
where service is offered to the customer class in which the queue is shortest.
The practical significance of this policy has been recognized \cite{nasser2005, urvoy2008, Giovanna2010, bonald2011, benameur2013, Guillemin2013, Guillemin2014, Guillemin2017}, and
analytic results have been obtained \cite{Guillemin2013, Guillemin2014, Guillemin2017, Giovanna2010}.
Briefly, our probabilistic assumptions are that both arrival and potential
service processes are Poisson, which makes the model Markovian,
and that arrival and service rates are class-dependent.
The diffusion scale behavior of the model in heavy traffic has not been
studied before. The main result of this paper addresses the $N$-dimensional
nominal workload (a term adopted from \cite{PKH}, expressing conditional expectation
of workload given the state)
and queue length processes, where $N$ denotes the number of classes.
It asserts that, under a critical load condition,
the diffusion scale versions of both these processes converge to processes
living in the set $\calS_0$, which consists of the union of
the $N$ coordinate axes in $\R_+^N$. Specifically, the rescaled nominal workload
converges to a {\it Walsh Brownian motion} (WBM) on $\calS_0$,
and the rescaled queue length converges to a certain diagonal transformation
of the same process.

WBM was introduced by Walsh in \cite{walsh1978} as a planar diffusion
that has a singular behavior at the origin.
Away from the origin it evolves as a one-dimensional Brownian motion (BM)
along a ray connecting its position to the origin, and its excursions
into rays emanating from the origin follow a fixed angular distribution.
Some early results on this process, including its special case referred to as
{\it skew BM}, where the state space consists of exactly two rays,
are \cite{Harrison1981, rogers1983ito, baxter1984equivalence, varopoulos1985long, salisbury1986construction, bar-pit-yor}.
Intriguing aspects related to the natural filtration of this process
were addressed in \cite{Tsirelson1997}.
Recently, vast extensions of this model have been proposed
and thoroughly studied. The reader is referred to \cite{ichiba2015stochastic} and the references therein
for this development.

In the terminology of Harrison \cite{har-bal}, an {\it unconventional}
limit theorem for a queueing system in heavy traffic
is one for which the limit process is not given as a reflected Brownian motion (RBM).
Our result thus belongs to a family of unconventional heavy traffic limits,
starting from \cite{har-wil-96} and including the more recent \cite{kruk-11}
as well as several other results surveyed in \cite{wil-96} and \cite{kruk-11}.
Moreover, our heavy traffic result is nonstandard in
that it does not constitute an {\it invariance principle}.
That is, it is observed in simulations that
the limit law (specifically, the angular distribution of the limit WBM)
is not determined solely by the first two moments of the data.
The simulations also indicate that it is sensitive even to tie breaking rules.
A third nonstandard aspect of the result is that the proof method does not
provide an explicit expression or a characterization of the limit law.
Whereas the modulus is given as an RBM with specified drift and diffusion coefficients,
no information on the angular distribution is available from the proof.
In fact, it appears unlikely to the authors that an explicit expression
can be attained except under some special symmetry.

Some further details on the policy are as follows.
In the literature, there are two variants,
distinguished by the interpretation given to the selection
of jobs from the shortest queue: that may refer to the one
having least nominal workload or the one having least number of jobs.
We adopt here the convention of \cite{Guillemin2013, Guillemin2014, Guillemin2017} and work with the former.
However, for all other purposes, the term {\it queue length} refers in this paper
to job count.
Next, the service rule is assumed to follow a preemptive priority.
Finally, the tie breaking rule is a part of the model description.
We allow for a rather general choice, by assuming that
when the collection, $\calK$, of classes having shortest queue consists of more
than one class, the server's effort is split according to a specified probability
measure $p^\calK$ supported on the set $\calK$.

Under static priority it is well known since Whitt's result \cite{whitt71}
that in heavy traffic, the queue which has least priority
is always {\it dominant}, where this term means that nearly all the workload
in the system is contained in this queue.
Under SSQ, heuristically, one may imagine that very soon after each time
the system (nearly) empties, a competition takes place among the queues, where
the one that loses ends up with most workload and consequently least priority.
Thus it is reasonable, in view of the aforementioned result on static priority,
to expect that the losing class
actually becomes dominant and remains so until again the system
becomes empty (or nearly empty). This heuristic suggests, moreover,
that the choice of the class to become dominant during an excursion
(the outcome of the competition, one may say) is random, and
is highly sensitive to
the dynamics of the Markov process as the queues just start to build.
The result of this paper reveals an asymptotic behavior
with exactly these elements.
One of the most significant and least obvious aspects of it
is that the probabilities of each class becoming dominant
starting from an empty system do converge in the scaling limit.
Indeed, their limit is given by the WBM's angular distribution.

An important feature of SSQ is that when two streams of arrivals have similar
first order characteristics but one is more variable than the other,
or has greater tendency to exhibit bursts,
the policy tends to prioritize the former over the latter.
This is due to the fact that a burst of traffic is likely
to cause a long queue, resulting in lower priority.
For this reason, SSQ has been referred to in the literature as {\it `implicit
service differentiation'} \cite{urvoy2008,bonald2011,Giovanna2010}
and {\it `self prioritization'} \cite{bonald2011}.
Quoting from \cite{Guillemin2014},
{\it ``...priority is thus implicitly given to smooth flows over data traffic...
sending packets in bursts''.}
The policy has gained interest in technological uses, specifically in the context of
packet scheduling \cite{urvoy2008, Giovanna2010, bonald2011, benameur2013, Guillemin2013, Guillemin2014, Guillemin2017}.
For example, in \cite{urvoy2008,bonald2011,benameur2013}
SSQ (referred to there as {\it shortest queue first})
is compared with
{\it first in first out} and {\it stochastic fairness queueing},
via experimental tests, and is argued to be the best candidate solution for
{\it quality of service} on ADSL internet access in various tests
(web browsing, file download, peer-to-peer file sharing,
VoIP and video calls, audio streaming, and video streaming).
It is also found experimentally that the policy
prioritizes TCP acknowledgment and delay- and loss-sensitive
applications (voice, audio and video streaming), which leads
to lower loss counts and delays.
For further advantages and additional uses of this policy see
\cite{Giovanna2010} and the references therein, as well as \cite{nasser2005}.

The policy has been theoretically analyzed in several papers.
Guillemin and Simonian \cite{Guillemin2014}
study the case of two buffers with Poisson arrivals and general service
time distributions, establishing
functional equations for the Laplace transform of the workload
processes at stationarity. They also specialize to the symmetric,
exponential service time case,
where they are able to derive empty queue probabilities
and tail behavior for the distribution of the workload.
In \cite{Guillemin2017} the authors study the same features in the asymmetric case,
again for $N=2$ at stationarity, where service times are exponentially distributed.
The paper \cite{Giovanna2010} studies
instantaneous throughput and buffer occupancy of $N\ge2$
long-lived TCP sources, using a deterministic fluid model,
under three per-flow scheduling disciplines:
fair queuing, longest queue first, and shortest queue first,
assuming longest queue drop buffer management.
They obtain closed form expressions for the stationary
throughput and the buffer occupancy.

We now make some comments about the proof.
To this end we introduce $\hat X^r(t)$, $t\in[0,\iy)$,
that are $[0,\iy)^N$-valued processes indexed by
the scaling parameter $r\in[1,\iy)$. The component $\hat X^r_i(t)$
represents the nominal workload in buffer $i$ at time $t$, rescaled diffusively;
the precise definition appears in \S \ref{sec2}.
We start by treating the rescaled total nominal workload,
$\sum_i\hat X^r_i(t)$, and recall the well-known
fact that it converges to an RBM under any
work conserving policy, to which SSQ is no exception.
This result is required in a slightly extended form,
stated in Lemma \ref{lem1}, which asserts that
convergence holds uniformly with respect to initial conditions.
The remainder of the proof has three main ingredients.
The first is concerned with showing that $\hat X^r$ resides
close to $\calS_0$ as $r$ gets large.
The aforementioned term `dominant queue'
is treated mathematically by considering tubes of width $\eps>0$
about each of the $N$ positive coordinate axes.
In terms of these tubes, queue $i$ is dominant at time $t$ if
$\hat X^r_i(t)$ resides in an $\eps$-tube about axis $i$,
for arbitrarily small $\eps$ and large $r$.
Thus the first main ingredient of the proof is to
show that the probability of exiting the collection of $N$ tubes
tends to zero as $r\to\iy$. This is the content of Lemma \ref{lem-b}(i).
Note that this element, along with the weak convergence of the total nominal workload to
an RBM, immediately provides the convergence of the modulus process
to the same RBM.

The second main ingredient is concerned with the angular behavior.
It is to show that the entrance law into tubes converges in the scaling limit.
We consider first a special case of the model,
that we call the {\it homogeneous} case, in which the transition intensities
of the underlying Markov process corresponding to $r>1$ are
rescaled version of those for $r=1$.
This trick buys us the ability to transform the double limit problem
of entrance law into $\eps$-tubes (involving $\eps\to0$ and $r\to\iy$)
to a single limit (involving $r\to\iy$ only).
The existence of a limit of the entrance law
is shown by arguing that, starting at the origin, the probabilities
of entering $r^{-\kap_0}$-tubes
form a Cauchy sequence, where $\kap_0>0$ is a suitable constant.
The tools used to establish this argument are
the martingale property of the total nominal workload (that also owes
to homogeneity), and a strengthening of Lemma \ref{lem-b}(i)
which improves the $o(1)$ exit probability estimates to polynomial
estimates.
Relying on the homogeneous case, the general case is then treated
by means of a change of measure.
The homogeneous case is stated in Lemma \ref{lem5}.
The double limit assertion is stated as Proposition \ref{lem3},
and the reduced version in the form of a single limit is given in \eqref{15}.
The polynomial exit probability measure is proved by means of
construction of a Lyapunov function for the distance of the state
from $\calS_0$, that may be interpreted as the nominal workload included
in all but the dominant class. This tool is stated in Lemma \ref{lem2}.
Finally, the change of measure argument is provided within the proof of
Proposition \ref{lem3} in \S \ref{sec35}.

The third main ingredient is the asymptotic independence of modulus and angle.
This relies, first and foremost, on the second ingredient alluded to above,
as well as on strong Markovity of the prelimit process and some estimates
on the heat kernel associated with RBM on the half-line.
This asymptotic independence property is stated in \eqref{22}.
These ingredients are finally combined in the proof of the main result,
building on the characterization of WBM via its semigroup \cite{bar-pit-yor},
and using crucially strong Markovity of the prelimit.

Some earlier results on the convergence of discrete processes to
WBM appear in \cite{Harrison1981} and \cite{Hajri2012}.
The paper \cite{Harrison1981} studies the case of a skew BM. The convergence
result included within this paper addresses a suitably defined
random walk on the integers observed
at the diffusion scale, and establishes its weak convergence to a skew BM.
The focus of \cite{Hajri2012} is the stochastic flow associated to WBM,
and for this model, discrete approximations to the flow are obtained.
In both these references,
the pre-limit processes already live in a collection of $N$
rays ($N=2$ in the former, $N\ge2$, finite, in the latter),
forming a symmetric random walk everywhere on the state space except
at the origin.
Consequently, the three main issues alluded to above in the description of our proof
(estimates on exiting tubes, existence of a limit for the entrance probability
into tubes, asymptotic independence) are all trivial in the cases studied
in \cite{Harrison1981} and \cite{Hajri2012}.

A general method was introduced in \cite{lamsim} for obtaining convergence of regenerative processes
from a certain notion of convergence of their excursions.
The regenerative processes we treat do fall into the category of those addressed
in \cite{lamsim}. However, in the setting considered here,
proving the convergence of excursions
amounts, roughly speaking, to establishing the three ingredients alluded to above,
and so it seems that as far as our result is concerned, this method does not provide
a significant shortcut.

The paper is organized as follows.
\S \ref{sec2} presents the model and the main result.
\S \ref{sec3} is devoted to the proof. First, in \S \ref{sec31+},
the result is proved based on Lemma \ref{lem1}, Lemma \ref{lem-b}
and Proposition \ref{lem3}, stated in the beginning of the section.
The convergence of the total nominal workload to an RBM is proved in \S \ref{sec32}.
\S \ref{sec33} provides
estimates on probabilities to exit the tubes.
\S \ref{sec34} and \S \ref{sec35} establish the limit result regarding
the angular distribution, dealing with the homogeneous case and
the general case, respectively.
Finally, some concluding remarks are included in \S \ref{sec4}.

\subsection*{Notation}
For $x,y\in\R^N$ ($N$ a positive integer), let $x\cdot y$
and $\|x\|$ denote the usual scalar product and $\ell_2$ norm, respectively.
Denote $[N]=\{1,2,\ldots,N\}$ and let $\{e_i:i\in[N]\}$
denote the standard basis in $\R^N$. Let $\one$ denote the
$N$-dimensional vector whose all entries equal 1.
For $x\in\R^N$ and $A\subset\R^N$, let $\dist(x,A)=\inf\{\|x-y\|:y\in A\}$.
Let $B(x,r)=\{y\in\R^N:\|y-x\|\le r\}$ denote the closed ball.
Denote $\R_+=[0,\iy)$.
For $f:\R_+\to\R^N$ and $T\in\R_+$, let $\|f\|_T=\sup_{t\in[0,T]}\|f(t)\|$,
and, for $\theta>0$,
$w_T(f,\theta)=\sup_{0\le s<u\le s+\theta\le T}\|f(u)-f(s)\|$.
For a Polish space $E$, let $\calC_E[0,T]$ and $\calD_E[0,T]$
denote the set of continuous and, respectively, c\`adl\`ag functions $[0,T]\to E$.
Let $\calC_E[0,\iy)$ and $\calD_E[0,\iy)$ denote the respective
sets of functions $[0,\iy)\to E$.
Endow $\calD_E[0,\iy)$ with the Skorohod $J_1$ topology.
A sequence of processes $\{X_n\}_n$ with sample paths in
$\calD_E[0,\iy)$ is said to be {\it $\calC$-tight} if it is tight and every subsequential
limit has, with probability 1, sample paths in $\calC_E[0,\iy)$. Write $X_n\To X$ for convergence in
law.
Let $\calC_0(E)$ denote the set of continuous, compactly supported functions
on $E$.
For $b\in\R$ and $\sig\in(0,\iy)$,
a {\it $(b,\sig)$-BM starting from $x\in\R$} is
a $1$-dimensional BM having drift $b$, infinitesimal covariance $\sig^2$
and initial condition $x$. A $(b,\sig)$-RBM starting from $x\in\R_+$
is an RBM in $\R_+$ with reflection at zero,
with the corresponding parameters and initial condition $x$.
Denote by $\calM_1$ the collection of $N$-dimensional
probability vectors, namely $\calM_1=\{x\in\R_+^N:\sum_ix_i=1\}$.
Throughout, we use the letter $c$ to denote a positive deterministic constant
whose value may change from one appearance to another.

\section{Setting and result}\label{sec2}

\subsection{Serve-the-shortest-queue in heavy traffic}\label{sec21}

Consider a sequence of queueing models indexed by $r\in[1,\iy)$, defined on a probability
space $(\Om,\calF,\bP)$.
A server operates to serve customers of $N\ge2$ classes.
Each customer class has a dedicated buffer with infinite room.
Upon arrival, a class-$i$ customer is queued in buffer $i\in[N]$.
The process representing the number of customers in buffer $i$
is called the {\it $i$th queue length}
and is denoted by $Q^r=(Q^r_1,\ldots,Q^r_N)$.
The $\Z_+^N$-valued random variable (RV) $Q^r(0)=(Q^r_1(0),\ldots,Q^r_N(0))$ is referred to as
the {\it initial queue length}.
The arrivals are Poissonian and the service times are exponential.
To model these, let $\{A^r_i\}_{i\in[N]},\{S^r_i\}_{i\in[N]}$
be a collection of $2N$ mutually independent Poisson processes,
with right-continuous sample paths, independent of the initial queue length,
where $A^r_i$ (resp., $S^r_i$) has rate $\la^r_i$ (resp., $\mu^r_i$).
The processes $A^r_i$ and $S^r_i$ represent the arrival and potential service
processes for class $i$, respectively. More precisely,
$A^r_i(t)$ is the number of class-$i$ customers to arrive
(to buffer $i$) until time $t$, and $S^r_i(t)$ gives the number of class-$i$
service completions by the time the server has dedicated $t$ units of time
to class-$i$ customers.

The process $X^r=(X^r_1,\ldots,X^r_N)$ defined by
\begin{align}\label{asaf1}
X^r_i=(\mu^r_i)^{-1}Q^r_i
\end{align}
is referred to as the {\it nominal workload process}.
This term, borrowed from \cite{PKH}, expresses the fact that $X^r_i(t)$ represents
the conditional expectation of the time it takes to serve
the $Q^r_i(t)$ customers present in buffer $i$ at time $t$,
conditioned on $Q^r_i(t)$
(assuming that the server works exclusively on this class).

Within each class, only one customer may be served at a time (and for concreteness,
we may assume it is the oldest one present in the system),
although service effort is sometimes split among classes (see below).
The priority rule among classes is to always serve the shortest queue
as measured in terms of nominal workload.
To make this statement precise some additional notation is required.
We say that {\it buffer $i$ contains the shortest queue}
at time $t$, if
\[
0<X^r_i(t)=\min\{X^r_j(t):X^r_j(t)>0,\, j\in[N]\}.
\]
When there is exactly one buffer containing
the shortest queue, the server serves it at full capacity (thus, service
is preemptive). When there is more than one such buffer, the server's effort
is split among the buffers containing the shortest queue
according to predetermined fractions in a head-of-the-line form. To model these fractions, it is assumed
that for any $\emptyset\ne\calK\subseteq[N]$ we are given a vector
$p^\calK\in\R_+^\calK$, such that $\sum_{i\in\calK}p^\calK_i=1$.
When the collection of shortest queues is $\calK$, the fraction of effort
dedicated to class $i$ is given by $p^\calK_i$. If we denote by $T^r_i(t)$
the total effort dedicated to class $i$ by time $t$ (measured in
units of time), then it is given by
\begin{equation}\label{37}
T^r_i(t)=\int_0^tp_i^{\calK(X^r(s))}ds,
\end{equation}
where, for $x\in\R_+^N$, we denote
\begin{equation}\label{39}
\calK(x)=\{i\in[N]:0<x_i\le x_j \text{ for all } j\in[N]\}.
\end{equation}
The departure process $D^r=(D^r_1,\ldots,D^r_N)$
consists of $N$ counting processes, where for each $i$, $D^r_i$ gives
the number of class-$i$ job completions. It thus satisfies
\begin{equation}\label{38-}
D^r_i(t)=S^r_i(T^r_i(t)).
\end{equation}
Clearly, $Q^r$ satisfies the balance equation
\begin{equation}
  \label{38}
  Q^r_i(t)=Q^r_i(0)+A^r_i(t)-D^r_i(t).
\end{equation}

This completes the description of the model.
Note that according to this description, the queue length process $Q^r$
is a Markov process on $\Z_+^N$, whereas $X^r$ is a Markov process on
\begin{equation}\label{28-}
\calS^r_{\rm u}=\frac{1}{\mu^r_1}\Z_+\times \cdots\times \frac{1}{\mu^r_N}\Z_+,
\end{equation}
(where `u' is mnemonic for {\it unscaled}).
Thus an alternative, concise description of the model is via the generator
of the process $X^r$, denoted by $\calL^r_{\rm u}$. It is given by
\begin{equation}
  \label{10}
  \calL^r_{\rm u}f(x)=\sum_{i\in[N]}\la^r_i\Big(f\Big(x+\frac{e_i}{\mu^r_i}\Big)-f(x)\Big)
  +\sum_{i\in\calK(x)}p^{\calK(x)}_i\mu^r_i\Big(f\Big(x-\frac{e_i}{\mu^r_i}\Big)-f(x)\Big),
\end{equation}
for any bounded $f:\calS^r_{\rm u}\to\R$.
Note that $\calK(0)=\emptyset$ and that, by the assumptions on $p^\calK$,
$p^\calK_i=1$ whenever $\calK$ consists of the singleton $\{i\}$, $i\in[N]$.

The parameters $\la^r_i$ and $\mu^r_i$ are assumed to scale like $r^2$.
The precise assumption is
that there exist constants $\la_i,\mu_i\in(0,\iy)$ and $\hat\la_i,\hat\mu_i\in\R$,
such that for $i\in[N]$, as $r\to\iy$,
\begin{align}\label{204}
\begin{split}
r^{-1}(\la^r_i-r^2\la_i)\to\hat\la_i,\\
r^{-1}(\mu^r_i-r^2\mu_i)\to\hat\mu_i.
\end{split}
\end{align}
The system is assumed to be critically loaded in the sense that the overall traffic intensity
equals 1. This is expressed as a condition on the first order parameters as follows,
\begin{align}\label{11}
\sum_{i\in[N]}\frac{\la_i}{\mu_i}=1.
\end{align}

Our main result regards rescaled versions of the nominal workload and queue length processes,
defined as 
\begin{equation}\label{14}
\hat X^r(t)=rX^r(t),
\qquad
\hat Q^r(t)=r^{-1}Q^r(t),\qquad t\in\R_+.
\end{equation}
Both these processes are obtained from $Q^r$ via
invertible transformations,
and are therefore Markov processes on discrete spaces.
The one to which most of the analysis is devoted in this paper is
$\hat X^r$. Recalling \eqref{asaf1}, it follows that $\hat X^r$ is a Markov process with state space
\begin{equation}\label{28}
\calS^r=\frac{r}{\mu^r_1}\Z_+\times \cdots\times \frac{r}{\mu^r_N}\Z_+.
\end{equation}
Specifically, the jump rates of both $X^r$ and $\hat Q^r$
are of order $r^2$ and their jump sizes are of order $r^{-1}$, confirming that \eqref{14}
gives the usual heavy traffic scaling.

\subsection{Walsh Brownian motion}\label{sec22}

In \cite{walsh1978}, Walsh introduced a diffusion process in the plane
that can informally be described as follows.
Let $\xi(t)=(\rho(t),\theta(t))$, $t\in\R_+$, be the representation of the process
in polar coordinates. Then the radial part $\rho(t)$
is an RBM, and on each excursion of $\xi(t)$ away
from the origin, the angular part $\theta(t)$ remains fixed.
Moreover, the constant value which $\theta(t)$ takes on each such excursion
has a fixed distribution, independent for the different excursions.
The precise definition that we shall work with is the one given
by Barlow, Pitman and Yor \cite{bar-pit-yor}, via its semigroup.
However, rather than working with a planar diffusion we work with what
is more natural for our purposes, namely a process in $\calS:=\R_+^N$.
Also, it is not necessary for our purposes
to consider general angular measures, and so the presentation
below is restricted to
angular measures supported on the $N$ vectors $\{e_i:i\in[N]\}$.

Let $b\in\R$, $\sig\in(0,\iy)$, and $q\in\calM_1$ be given.
Let $\Pi_t^+$, $t\in\R_+$ and $\Pi_t^0$, $t\in\R_+$ denote
the semigroups of a $(b,\sig)$-RBM and a $(b,\sig)$-BM killed at 0, respectively.
That is, for $f\in \calC_0(\R_+)$,
\[
\Pi^+_tf(x)=\Erbm_x[f(\rho(t))],\qquad \Pi_t^0f(x)=\Erbm_x[f(\rho(t))\one_{\{t<\zeta\}}],\qquad x\in\R_+,
\]
where $\rho(t)$ is a $(b,\sig)$-RBM and $\zeta$ denotes its hitting time at zero,
and, throughout the paper,
$\PPrbm_x$ (resp., $\Erbm_x$) denotes the law of $\rho$ with $\rho(0)=x$
(resp., the corresponding expectation).
Let $S^k$ denote the $k$-sphere. Use polar coordinates $(\rho,\theta)\in\R_+\times S^{N-1}$
to denote members $x\in\calS$
by setting $\rho=\|x\|$ and $\theta=x/\|x\|$ when $x\ne0$, $\theta=e_1$ when $x=0$.
The semigroup $\PI_t$ of a $(b,\sig,q)$-WBM is defined as follows.
For $f\in \calC_0(\calS)$, $\PI_t$ acts on $f$ as
\begin{equation}\label{9}
  \begin{split}
  & \PI_tf(0,\theta)=\Pi^+_t\bar f(0),
  \\
  & \PI_tf(\rho,\theta)=\Pi^+_t\bar f(\rho)+\Pi^0_t(f_\theta-\bar f)(\rho),
  \end{split}
\end{equation}
where we denote
\begin{equation}\label{9a}
\begin{split}
  & \bar f(\rho)=\sum_{i\in[N]}q_if(\rho,e_i),\qquad \rho\ge0,\\
  & f_\theta(\rho)=f(\rho,\theta),\qquad\qquad\ \rho\ge0,\ \theta\in S^{N-1}.
\end{split}
\end{equation}
It is shown in \cite{bar-pit-yor} that $\PI_t$ is a Feller semigroup on $\calC_0(\calS)$
and that there exists a strong Markov process $\{\xi(t)\}$
with state space $\calS$ and semigroup $\PI_t$, that has a.s.-continuous sample paths.
Moreover, this process has the properties alluded to above. More precisely,
when written in polar coordinates as $\xi(t)=(\rho(t),\theta(t))$,
the radial part $\rho(t)$ is a $(b,\sig)$-RBM and the values that
the angular part $\theta(t)$ takes are constant on the interval $[0,\zeta]$
(where the constant is determined by the initial condition $\theta_0$)
as well as on each excursion away from zero.
These constant values on the excursions away from zero are mutually independent
with common distribution $\sum_{i\in[N]}q_i\del_{e_i}(dx)$, where $\delta_{e_i}$ is the Dirac measure at $e_i$.
In this paper we are interested in the case where the initial condition
is supported on $\calS_0:=\cup_{i\in[N]}\{xe_i:x\in\R_+\}$.
Note that in this case, $\xi(t)$ takes values in $\calS_0$ for all $t$.

Throughout, let $\PP_x^{\rm wbm}$ and $\E_x^{\rm wbm}$ denote the law
of $\xi$ for $\xi(0)=x$, and respective expectation.
Then relations \eqref{9} can be expressed,
for $x=(\rho_0,\theta_0)$, as
\begin{equation}\label{8}
\E_x^{\rm wbm}[f(\xi(t))]=\E_{\rho_0}[f(\rho(t)\theta_0)\one_{\{t<\zeta\}}]
+\sum_{i\in[N]}q_i\E_0[f(\rho(t)e_i)\one_{\{t\ge\zeta\}}].
\end{equation}

\subsection{Main result}\label{sec23}

The linear relation between $X^r$ and $Q^r$, the convergence $r^{-2}\mu^r_i\to\mu_i$
that follows from \eqref{204}, and the rescaling defined in \eqref{14} imply
an asymptotic relation between $\hat X^r$ and $\hat Q^r$ which one can express
in terms of the $N\times N$ matrix $\hat M=\diag(\mu_i)_{i\in[N]}$.
For example, the statement $\hat X^r(0)\To\xi(0)$ is equivalent to
the statement $\hat Q^r(0) \To \hat M\xi(0)$, as $r\to\iy$,
where, throughout, the symbol $\To$ denotes convergence in law under $\bP$. Denote
\[
b=\sum_{i\in[N]}\frac{1}{\mu_i}\Big(\hat\la_i-\frac{\la_i}{\mu_i}\hat\mu_i\Big),
\qquad
\sig^2=2\sum_{i\in[N]}\frac{\la_i}{\mu_i^2}.
\]

\begin{theorem}
  \label{th1}
There exists $q\in\calM_1$ such that, if
$\{\xi(t)\}$ is a $(b,\sig,q)$-WBM with initial distribution supported
on $\calS_0$ and $\hat X^r(0)\To\xi(0)$
then $\hat X^r\To\xi$ and $\hat Q^r\To\hat M\xi$, as $r\to\iy$.
\end{theorem}

\begin{remark}
(a) Whereas the coefficients $b$ and $\sig$ of the process $\xi$
are given explicitly, our approach does not provide
a construction or any explicit information of the angular distribution $q$.
However, this much can be said: $q$ does not depend on the second order
parameters $(\hat\la_i,\hat\mu_i)$
(where we use standard terminology by which $(\la_i,\mu_i)$ and $(\hat\la_i,\hat\mu_i)$ are called
first and second order parameters, respectively, due to the fact that in most conventional queueing models,
LLN limits depend only on the former, whereas CLT limits are also affected by the latter).
This statement is a direct consequence of our results of \S\ref{sec35}.
\\
(b) Initial conditions which are not asymptotically concentrated
on $\calS_0$ are excluded from our treatment.
For such initial conditions the asymptotic behavior
is expected to follow a jump to $\calS_0$ at time zero,
and then proceed as a WBM.
However, the position to which the process jumps is dictated by
properties finer than the limiting initial distribution,
to the extent that the limit does not exist in general.
For example, for $N=2$,
a sequence of initial conditions may converge to a point
on the diagonal in such a way that $\hat Q^r_1(0)>\hat Q^r_2(0)+\eps_r$.
It is not hard to see that,
due to even a small advantage $\eps_r>0$ to $\hat Q^r_2(0)$,
the limiting process will initially jump to a point on the $e_1$ axis,
provided that $\eps_r$ tends to zero sufficiently slowly.
Interchanging the roles of $\hat Q^r_1(0)$ and $\hat Q^r_2(0)$
will result in a jump to the $e_2$ axis.
\end{remark}

\section{Proof of the main result}\label{sec3}

Below we present two central lemmas and one central proposition required to prove our main result.
The proof of the main result is presented next, in \S \ref{sec31+}.
The proofs of the lemmas and the proposition are then provided in \S \ref{sec32}--\ref{sec35}.

Some notation used throughout this section is as follows.
We use $\sum_i$ as shorthand notation for $\sum_{i\in[N]}$.
For $\ph\in\calD_\R[0,\iy)$, let $\Gam[\ph]=(\Gam_1[\ph],\Gam_2[\ph])$ be defined by
\begin{equation}\label{40}
(\Gam_1[\ph](t),\Gam_2[\ph](t))=\Big(\ph(t)-\inf_{s\le t}(\ph(s)\w0),\
-\inf_{s\le t}(\ph(s)\w0)\Big),\qquad t\in[0,\iy).
\end{equation}
The {\it Skorohod map} $\Gam$ just introduced
transforms a $(b,\sig)$-BM starting from $x\ge0$,
say, $W$, into a $(b,\sig)$-RBM starting from the same point, via
$R=\Gam_1[W]$. The process given by $\Gam_2[W]$ gives the corresponding boundary term.

Let $\hat R^r(t)=\one\cdot\hat X^r(t)$, $t\in\R_+$, and let $\rho$
be a $(b,\sig)$-RBM. In addition to the notation $\PP^{\rm wbm}_x$ and $\PP_x$ introduced
above, for each $r$ and $x\in\calS^r$,
we use $\PP^r_x$ and $\E^r_x$ for the law of the Markov process
$\hat X^r$ with $\hat X^r(0)=x$, and the respective expectation.
Moreover, for each $r$ and $x\in\calS^r$,
we use $\bP^r_x$ and $\bE^r_x$ for the law of the tuple
$(A^r,S^r,X^r)$ with $\hat X^r(0)=x$,
and the respective expectation.

Let $\calS_\eps=\{x\in\calS:\dist(x,\calS_0)<\eps\}$.
Finally, denote $\zeta^r=\inf\{t\ge0:\hat R^r(t)=0\}$, and for $\eps>0$,
\begin{equation}\label{35}
\tau^r(\eps)=\inf\{t\ge 0:\hat R^r(t)\ge\eps\}.
\end{equation}
Both $\zeta^r$ and $\tau^r(\eps)$ are easily seen to be a.s.\ finite.

\begin{lemma}
	\label{lem1}
i.
The process $\hat R^r$ is given as $\hat R^r=\Gam_1[\hat R^r(0)+B^r]$,
where $B^r$ decomposes as $B^r=\tilde B^r+E^r$.
For each $r$, $\tilde B^r$ (resp., $E^r$)
is measurable w.r.t.~$\sig\{A^r(t),S^r(t),t\in\R_+\}$
(resp., $\sig\{A^r(t),S^r(t), X^r(t), t\in\R_+\}$)
and $\tilde B^r\To B$, where $B$ is a $(b,\sig)$-BM starting from zero,
whereas
\begin{align}\label{cohen1}
\lim_{v\down0}\limsup_{r\to\iy}\sup_{x\in\calS^r}\bP^r_x(\|E^r\|_T>v)=0.
\end{align}
As a consequence, if $\hat R^r(0)\To \rho(0)$, then $\hat R^r\To\rho$.
\\
	ii. For any $t_0>0$,
	\[
	\lim_{v\down0}\liminf_{r\to\iy}\inf_{x\in\calS^r:\one\cdot x<v}
	\PP^r_x(\zeta^r\le t_0)=1.
	\]
\end{lemma}

Throughout, let $\calU_0$ denote the class of functions
$u:[1,\iy)\to(0,\iy)$ for which $u(r)\to0$ as $r\to\iy$.

\begin{lemma}
  \label{lem-b}
The processes $\hat X^r$ are $\calC$-tight under $\bP$.
Moreover, let $\nu^r$ denote the distribution of $\hat X^r(0)$.
Then there exists $u\in\calU_0$ such that
for every $T>0$ one has the following.
\\
i.
$\PP^r_{\nu^r}(\hat X^r(t)\in\calS_{u(r)} \text{ for all } t\in[0,T])\to1$ as $r\to\iy$.
\\
ii.
$\inf\PP^r_x(\hat X^r(t)\in\calS_{u(r)} \text{ for all } t\in[\zeta^r\w T,T))\to1$ as $r\to\iy$,
where the infimum extends over $x\in\calS^r\cap K$, and $K\subset\calS$ is a given
compact set.
\\
iii. For $f\in\calC_0(\calS)$, $t\in\R_+$, $i\in[N]$, and $k\in(0,\iy)$,
\[
\lim_{\del\down0}
\limsup_{r\to\iy}
\sup_{y\in\calS^r,x\in[0,k]:\|y-xe_i\|<\del}
\Big|\E^r_y[f(\hat X^r(t))\one_{\{t<\zeta^r\}}]
-\Erbm_x[f(\rho(t)e_i)\one_{\{t<\zeta\}}]\Big|=0,
\]
where we recall that $\zeta=\inf\{t\ge 0:\rho(t)=0\}$.
\end{lemma}

\begin{proposition}
\label{lem3}
There exist $q\in\calM_1$ and $u\in\calU_0$ such that
\begin{align}\label{200}
\lim_{\eps\down0}\;\limsup_{r\to \iy}|\PP^r_0(\hat X^r(\tau^r(\eps))
\in B(\eps e_i,u(r))-q_i|=0,\qquad i\in[N].
\end{align}
\end{proposition}

\subsection{Proof of Theorem \ref{th1}}\label{sec31+}

Given $\eps>0$,
define a sequence of hitting times as
\begin{align}\label{29}
\notag
\zeta_0^r&=\inf\{t\ge0:\hat R^r(t)=0\},\\
\tau^r_m&=\inf\{t\ge\zeta^r_m:\hat R^r(t)\ge \eps\},
\qquad m=0,1,\ldots
\\
\notag
\zeta^r_{m+1}&=\inf\{t\ge\tau^r_m:\hat R^r(t)=0\},
\qquad m=0,1,\ldots
\end{align}
Let $N^r_t=\sup\{m:\tau^r_m\le t\}$.
When we need to emphasize the dependence on $\eps$ we write these RVs as $\zeta^r_m(\eps)$
and $\tau^r_m(\eps)$.

Let $(\xi(t))_{t\in\R_+}$ be a $(b,\sig,q)$-WBM and assume, without loss of generality,
that $\rho=\one\cdot\xi$.
For this process we define an analogous sequence of hitting times by
\begin{align*}
\zeta_0&=\inf\{t\ge0:\rho(t)=0\},\\
\tau_m&=\inf\{t\ge\zeta_m:\rho(t)\ge \eps\},
\qquad m=0,1,\ldots
\\
\zeta_{m+1}&=\inf\{t\ge\tau_m:\rho(t)=0\},
\qquad m=0,1,\ldots,
\end{align*}
and set $N_t=\sup\{m:\tau_m\le t\}$.

The weak convergence stated in Lemma \ref{lem1}(i) does not directly imply that
of the hitting times $\tau^r(\eps)$ of \eqref{35}
to $\tau(\eps):=\inf\{t\ge 0:\rho(t)\ge\eps\}$ when both $\hat R^r$ and $\rho$
start at zero.
However, this convergence is clearly valid, as can be seen by using
in addition the property of RBM that $\tau(\eps+\delta)\to\tau(\eps)$
in probability as $\delta\down0$.
Moreover, under $\bP$, it is assumed in Theorem \ref{th1}
that $\hat X^r(0)$ converges to $\xi(0)$ in distribution.
An inductive use of this fact yields a similar statement
for the stopping times $\{\tau^r_m\}_m$.
More precisely, for any fixed $m$, as $r\to\iy$, we have,
the following uniform convergence: for any compact set $K\subset\calS$ and a function $h\in\calC_0(\R_+)$,
\begin{equation}
  \label{34}
 \limsup_{\delta\down0}\limsup_{r\to\iy}
\sup_{\substack{x\in K\cap\calS_0}}
\sup_{\substack{y\in\calS^r \\ \|x-y\|<\delta}}
\left|\E^r_y[h(\tau^r_m)]-\Erbm_{\|x\|}[h(\tau_m)]\right|=0.
 \end{equation}

The proof of the main result
is based on finite-dimensional convergence and $\calC$-tightness.
The key ingredient is showing that for any compact set $K\subset \calS$,  $t\in\cal\R_+$,
and a function $f\in\calC_0(\calS)$,
\begin{align}\label{ucon}
\limsup_{\delta\down0}\limsup_{r\to\iy}
\sup_{\substack{x\in K\cap\calS_0}}
\sup_{\substack{y\in\calS^r \\ \|x-y\|<\delta}}
\left|\E^r_y[f(\hat X^r(t))]-\E^{\rm wbm}_x[f(\xi(t))]\right|=0.
\end{align}
Before proving this statement we show, adapting the proof of
Theorem 4.2.5 of \cite{ethkur}, that it implies the convergence of
$\hat X^r$ to $\xi$ for finite-dimensional marginals.
That is, for every $x\in\calS_0$, $\{x^r\}_r$, $x^r\in\calS^r$ that converges to $x$,
$m\ge 1$, $0\le t_1<\cdots<t_m$,
and functions $h_1,\ldots,h_m\in\calC_0(\calS)$, one has
\begin{align}\label{300}
\lim_{r\to\iy}\E^r_{x^r}[h_1(\hat X^r(t_1))\cdots h_m(\hat X^r(t_m))]=
\E^{\text{wbm}}_{x}[h_1(\xi(t_1))\cdots h_m(\xi(t_m))].
\end{align}
We argue by induction over $m$.
The base case follows from \eqref{ucon}.
Next, assume that \eqref{300} holds for $m$.
Denote by $\PI^r_t$ the semigroup corresponding to $\{\hat X^r(t)\}$.
Then by Lemma \ref{lem-b}(i), there exists $u\in\calU_0$, such that
\begin{align}\notag
&\E^r_{x^r}[h_1(\hat X^r(t_1))\cdots h_m(\hat X^r(t_m))\cdot h_{m+1}(\hat X^r(t_{m+1}))]\\\notag
&\quad=\E^r_{x^r}[h_1(\hat X^r(t_1))\cdots h_m(\hat X^r(t_m))\cdot \Pi^r_{t_{m+1}-t_m}h_{m+1}(\hat X^r(t_m))]\\\notag
&\quad=\E^r_{x^r}[h_1(\hat X^r(t_1))\cdots h_m(\hat X^r(t_m))\cdot \Pi^r_{t_{m+1}-t_m}h_{m+1}(\hat X^r(t_m))\one_{\{\hat X^r(t_m)\in\calS_{u(r)}\}}]+o_r(1),
\end{align}
where here $o_r(1)$ denotes a generic function of $r$ that vanishes as $r\to\iy$.
From \eqref{ucon} and the Feller property of $\PI_t$ proved in \cite{bar-pit-yor},
it follows that for $h\in\calC_0(\calS)$,
$\sup_{x\in\calS^r\cap\calS_{u(r)}}|\PI^r_th(x)-\PI_th(x)|\to0$ as $r\to\iy$.
It follows that the expression in the above display equals
\begin{align}\notag
&\E^r_{x^r}[h_1(\hat X^r(t_1))\cdots h_m(\hat X^r(t_m))\cdot \Pi_{t_{m+1}-t_m}h_{m+1}(\hat X^r(t_m))\one_{\{\hat X^r(t_m)\in\calS_{u(r)}\}}]+o_r(1)\\\label{301}
&\quad=\E^r_{x^r}[h_1(\hat X^r(t_1))\cdots h_m(\hat X^r(t_m))\cdot \Pi_{t_{m+1}-t_m}h_{m+1}(\hat X^r(t_m))]+o_r(1).
\end{align}
By the induction hypothesis, the above expression converges to
\begin{align*}
&\E^{\text{wbm}}_x[h_1(\xi(t_1))\cdots h_m(\xi(t_m))\PI_{t_{m+1}-t_m}h_{m+1}(\xi(t_m))]\\
&=
\E^{\text{wbm}}_x[h_1(\xi(t_1))\cdots h_{m+1}(\xi(t_{m+1}))].
\end{align*}
This establishes \eqref{300}.
In view of the $\calC$-tightness of $\hat X^r$ stated in Lemma \ref{lem-b},
this gives the main result $\hat X^r\To\xi$.

The rest of the proof is devoted to showing that \eqref{ucon} holds.
Fix $f\in\calC_0(\calS)$.
It suffices to prove the result for $f(0)=0$. Moreover, arguing by approximation,
we may, and will assume that $f$ is constant on a ball about the origin.
Thus, there exists $\eps>0$
for which $f(x)$ vanishes for all $x$ with $\one\cdot x\le \eps$.
We fix such $\eps$, and let $\tau^r_m=\tau^r_m(\eps)$,
and similarly let $\zeta^r_m$, $\tau_m$ and $\zeta_m$
be defined in terms of the same value of $\eps$.

Fix $t>0$ and a compact set $K\subset \calS_0$. For $u\in\calU_0$,
we will be concerned with $x\in K$ and $y^r\in\calS^r\cap\calS_{u(r)}$ such that
$\|x-y^r\|<u(r)$. We call such a pair $(x,(y^r)_r)$
a $u$-{\it admissible pair}. In what follows we denote $y=(y^r)$.
Since $t$ and $K$ are arbitrary, to prove \eqref{ucon},
it suffices to show that $\E^r_{y^r}[f(\hat X^r(t))]\to\E^{\rm wbm}_x[f(\xi(t))]$
uniformly over $u$-admissible pairs $(x,y)$, for an arbitrary $u\in\calU_0$.
Fix such a function $u$. Notice that the assertions in Lemma \ref{lem-b}(i), (ii), and Proposition \ref{lem3} are all monotone in $u$ in the sense that if they hold for some $u\in\calU_0$ then they also hold for a function that dominates $u$ and vanishes at infinity. Hence, without loss of generality, we may and will assume that Lemma \ref{lem-b}(i), (ii), and Proposition \ref{lem3} hold for the function $u$ that we have fixed.

On the intervals $[\zeta^r_m,\tau^r_m)$ one has $\hat R^r(t)\le \eps$.
As a consequence,
\[
\E^r_{y^r}[f(\hat X^r(t))]=F^{0,r}_{y}+F^r_{y},
\]
where
\[
F^{0,r}_{y}:=\E^r_{y^r}[f(\hat X^r(t))\one_{\{0\le t<\zeta^r_0\}}],
\qquad
F^r_{y}:=\sum_{m=0}^\iy \E^r_{y^r}[f(\hat X^r(t))\one_{\{\tau^r_m\le t<\zeta^r_{m+1}\}}].
\]
The above goal will be achieved once we show that, uniformly over $u$-admissible pairs,
\begin{align}\label{302}
F^{0,r}_{y}\to \Erbm_x[f(\rho(t)\theta)\one_{\{t<\zeta_0\}}]
\quad\text{and}\quad
F^r_{y}\to \Erbm_x[\bar f(\rho(t))\one_{\{t\ge\zeta_0\}}],
\end{align}
where we recall the definition of $\bar f$ from \eqref{9a}, that $\theta=x/\|x\|$ for $x\ne 0$ and $\theta=e_1$ for $x=0$.
Note that the first convergence is stated in Lemma \ref{lem-b}(iii).
Thus in what follows we focus on the term $F^r_y$. Denote $$\chi^r_m=\hat X^r({\tau^r_m}).$$
Recall that the jumps of the (unscaled) queue length process $Q^r$ are of size 1.
By the way the scaled nominal workload process $\hat X^r$ is defined, it follows
that all the jumps of this process are bounded by $cr^{-1}$, for some positive constant $c$.
As a consequence, one always has $\eps\le \|\chi^r_m\|\le \eps+cr^{-1}$. Denote
$B^r_i=B(\eps e_i,u(r))$.
It follows from Lemma \ref{lem-b}(ii)
that
\[
\PP^r_{y^r}(\text{for all } m\le N^r_t, \chi^r_m\in\cup_iB^r_i)\to1,
\]
uniformly over $u$-admissible pairs. As a result,
\[
F^r_y=\sum_i\sum_{m=0}^\iy \E^r_{y^r}[f(\hat X^r(t))\one_{\{\tau^r_m\le t<\zeta^r_{m+1}\}}
\one_{\{\chi^r_m\in B^r_i\}}]+o_r(1),
\]
where here and in what follows, $o_r(1)$ is a generic function of $r, x$ and $y$,
that converges to zero as $r\to\iy$, uniformly over $u$-admissible pairs
$(x,y)$.

Next we truncate the sum over $m$. The tail
$\sum_{m>M}$ can be estimated by $\|f\|_\iy \PP^r_{y^r}(N^r_t>M)$.
For a fixed initial condition $x$, the $\calC$-tightness of $\hat R^r$
gives the tightness of the RVs $N^r_t$. For arbitrary initial conditions $y^r$,
the strong Markovity reduces the same question to that of
tightness of $N^r_t$ when starting at $0$. Thus
\begin{equation}\label{12}
F^r_y=\sum_i\sum_{m=0}^M\E^r_{y^r}[f(\hat X^r(t))\one_{\{\tau^r_m\le t<\zeta^r_{m+1}\}}
\one_{\{\chi^r_m\in B^r_i\}}]+o_{M,r}(1),
\end{equation}
where here and in what follows, $o_{M,r}(1)$ refers to any function
$g$ of $(x,y,r,M)$
satisfying $\lim_{M\to\iy}\limsup_{r\to\iy}\sup_{(x,y)\text{ $u$-admissible}}|g(x,y,r,M)|=0$.

Our next step is to use the condition $\chi^r_m\in B^r_i$ included
in the $(i,m)$-th term in \eqref{12},
to approximate the expression $f(\hat X^r(t))$ therein by $f(\hat R^r(t)e_i)$.
Note carefully that it is possible for the process to move from $B^r_i$
to $B^r_j$, $j\ne i$, without exiting $\calS_{u(r)}$ or hitting the origin.
Thus we must argue that, given any distinct $i,j\in[N]$ and $t>0$,
\begin{align}
  \label{13}\notag
  &\PP^r_{y^r}(\text{there exists $m\in\{0,\ldots,M\}$ such that}\\
  &\qquad\tau^r_m\le t<\zeta^r_{m+1},
  \chi^r_m\in B^r_i,\hat R^r(t)>\eps,\|\hat X^r(t)-\hat R^r(t)e_j\|\le u(r))=o_r(1).
\end{align}
The proof of this statement, which we now give, is based on the fact that
in order for the process to behave as indicated in \eqref{13}
while remaining within $\calS_{u(r)}$, it must reach close to the origin.
Since we consider only finitely many $m$'s, it is sufficient
to show that for every fixed $t,m$, and $j\ne i$,
\begin{equation}
\label{ac002}
  \PP^r_{y^r}( \tau^r_m\le t<\zeta^r_{m+1},
  \chi^r_m\in B^r_i,\hat R^r(t)>\eps,\|\hat X^r(t)-\hat R^r(t)e_j\|\le u(r))=o_r(1).
\end{equation}
For every $r\in[1,\iy)$, $m\in\N$, and $j\in[N]$ define
\[
\pi^r_m=\pi^r_m[j]=\inf\{s\ge\tau^r_m:\hat R^r(s)>\eps,\|\hat X^r(s)-\hat R^r(s)e_j\|\le u(r)\}.
\]
The probability from \eqref{ac002} is bounded above by
$\PP^r_{y^r}( \tau^r_m\le t\le\pi^r_m<\zeta^r_{m+1},\chi^r_m\in B^r_i)$.
Under this event, if the process
does not leave $\calS_{u(r)}$ between times $\tau^r_m$ and $\zeta^r_m$,
then after time $\tau^r_m$, it must reach close to the origin
and not hit the origin prior to reaching a small neighborhood of $\eps e_j$.
Therefore, the LHS of \eqref{ac002} is bounded by
\begin{align}\notag
  &\PP^r_{y^r}( \exists s\in[0,t],\; \hat X^r(s)\notin\calS_{u(r)})+\PP^r_{y^r}(\phi^r_m<\pi^r_m<\zeta^r_{m+1}),
\end{align}
where $\phi^r_m=\inf\{s\ge\tau^r_m: \hat R^r(s)\le u(r)\sqrt{N} \}. $
From Lemma \ref{lem-b}(ii), the first term is $o_r(1)$.
For every $r$, let $\{\calF^r_s\}$ denote the filtration induced by
$\{\hat X^r(t)\}$.
Then for the second term, using the strong Markov property,
\begin{align*}
\PP^r_{y^r}(\phi^r_m<\pi^r_m<\zeta^r_{m+1})
&=
\E^r_{y^r}[\E^r_{y^r}[\one_{\{\phi^r_m<\pi^r_m<\zeta^r_{m+1}\}}\mid\calF^r_{\phi^r_m}]]\\
&=\E^r_{y^r}[\psi^r_j(\hat X^r(\phi^r_m))],
\end{align*}
where $\psi_j^r(z)=\E^r_z[\one_{\{\hat\pi^r< \zeta^r_0\}}]$
and $\hat\pi^r=\inf\{s\ge 0:\hat R^r(s)>\eps,\|\hat X^r(s)-\hat R^r(s)e_j\|\le u(r)\}$.
It remains to show that
$\lim_{r\to\iy}\sup_{\|z\|\le cu(r)}\psi^r_j(z)=0$,
for $c>0$ a constant.
Now, $\psi^r_j(z)\le
\E^r_z[\one_{\{\hat\tau^r< \zeta^r_0\}}]$, where
$\hat\tau^r=\inf\{s\ge 0:\hat R^r(s)\ge \eps\}$.
The last term goes to zero, uniformly in $\|z\|\le cu(r)$.
This shows \eqref{13}.

Equipped with \eqref{13}, we have from \eqref{12}
\[
F^r_y
=\sum_i\sum_{m=0}^M \E^r_{y^r}[f(\hat X^r(t))\one_{\{\tau^r_m\le t<\zeta^r_{m+1}\}}
\one_{\{\chi^r_m\in B^r_i\}}\one_{\{\|\hat X^r(t)-\hat R^r(t)e_i\|\le u(r)\}}]
+o_{M,r}(1).
\]
Thus, using the uniform continuity of $f$ and denoting $f_i(z)=f(ze_i)$, $z\in\R_+$,
\begin{align*}
F^r_y
&=\sum_i\sum_{m=0}^M \E^r_{y^r}[f_i(\hat R^r(t))\one_{\{\tau^r_m\le t<\zeta^r_{m+1}\}}
\one_{\{\chi^r_m\in B^r_i\}}\one_{\{\|\hat X^r(t)-\hat R^r(t)e_i\|\le u(r)\}}]+o_{M,r}(1)
\\
&=\sum_i\sum_{m=0}^M \E^r_{y^r}[f_i(\hat R^r(t))\one_{\{\tau^r_m\le t<\zeta^r_{m+1}\}}
\one_{\{\chi^r_m\in B^r_i\}}]+o_{M,r}(1),
\end{align*}
again using Lemma \ref{lem-b}.
The $(i,m)$-th term can be written as
\[
\E^r_{y^r}[\E^r_{y^r}[f_i(\hat R^r(t))
\one_{\{\tau^r_m\le t<\zeta^r_{m+1}\}}|\calF^r_{\tau^r_m}]
\one_{\{\chi^r_m\in B^r_i\}}].
\]
By strong Markovity, the conditional expectation
above can be written as $\ph^r_i(\tau^r_m,\chi^r_m)$, where
\begin{equation}\notag
\ph^r_i(s,z)=\E^r_{z}[f_i(\hat R^r(t-s)\one_{\{\zeta^r_0>t-s\}}]\one_{\{s\le t\}}.
\end{equation}
This gives
\[
F^r_y
=\sum_i\sum_{m=0}^M \E^r_{y^r}[\ph^r_i(\tau^r_m,\chi^r_m)\one_{\{\chi^r_m\in B^r_i\}}]+o_{M,r}(1).
\]
Define
$\ph_i(s)=\Erbm_\eps[f_i(\rho_{t-s})\one_{\{\zeta_0>t-s\}}]\one_{\{s\le t\}}$.
Then by Lemma \ref{lem-b}(iii)
one has
\(
\ph^r_i(s,z)=\ph_i(s)+o_r(1)
\) for $z\in B^r_i$.
Hence
\[
F^r_y
=\sum_i\sum_{m=0}^M \E^r_{y^r}[\ph_i(\tau^r_m)\one_{\{\chi^r_m\in B^r_i\}}]+o_{M,r}(1).
\]
It will be shown below that, for fixed $(i,m)$,
$\tau^r_m$ and $\chi^r_m$ are asymptotically independent, in the sense that
\begin{equation}
  \label{22}
\E^r_{y^r}[\ph_i(\tau^r_m)\one_{\{\chi^r_m\in B^r_i\}}]
=\E^r_{y^r}[\ph_i(\tau^r_m)]q_i+o_r(1).
\end{equation}
Hence
\[
F^r_y
=\sum_iq_i\sum_{m=0}^M \E^r_{y^r}[\ph_i(\tau^r_m)]+o_{M,r}(1).
\]
Using \eqref{34},
and a similar argument based on strong Markovity,
\begin{align*}
F^r_y
&=\sum_iq_i\sum_{m=0}^M \Erbm_{\|x\|}[\ph_i(\tau_m)]+o_{M,r}(1)
\\
&=\sum_iq_i\sum_{m=0}^M\E_{\|x\|}[f_i(\rho(t))\one_{\{\tau_m\le t<\zeta_{m+1}\}}]+o_{M,r}(1)
\\
&=\sum_iq_i\E_{\|x\|}[f_i(\rho(t))\one_{\{t\ge\zeta_0\}}]+o_{M,r}(1)
\\
&=\Erbm_{\|x\|}[\bar f(\rho(t))\one_{\{t\ge\zeta_0\}}]+o_{M,r}(1).
\end{align*}
Thus sending $r\to\iy$, then $M\to\iy$ gives the second
statement in \eqref{302}.

It remains to prove \eqref{22}.
Since $\ph_i$ is continuous on $[0,t]$ and the law of $\tau_m$ has no atoms,
it suffices to prove that for every $s\in[0,t]$ and $(i,m)$,
\begin{equation}
  \label{31}
  \PP^r_{y^r}(\tau^r_m\le s,\chi^r_m\in B^r_i)\to \PPrbm_{\|x\|}(\tau_m\le s)q_i,
\end{equation}
uniformly over $u$-admissible pairs $(x,y)$.
Toward showing \eqref{31}, we argue
that it suffices to establish this assertion for $y\equiv0$ and $m=0$.
Indeed,
\begin{align*}
\PP^r_{z}(\tau^r_m\le s,\chi^r_m\in B^r_i)
&=\E^r_{z}[\E^r_{z}[\tau^r_m\le s,\chi^r_m\in B^r_i|\calF^r_{\zeta^r_m}]]
\\&
=\E^r_{z}[\tilde\ph^r(\zeta^r_m)],
\end{align*}
where we used the fact that $\zeta^r_m<\tau^r_m$,
and $\hat X^r(\zeta^r_m)=0$, and denoted
\[
\tilde\ph^r(s)=\PP^r_0(\tau^r_0<t-s,\hat X^r(\tau^r_0)\in B^r_i).
\]
Similarly, $\PP^r_{z}(\tau^r_m\le s)=\E^r_{z}[\ph^{*,r}(\zeta^r_m)]$,
$\ph^{*,r}(s)=\PP^r_0(\tau^r_0<t-s)\to \PPrbm_0(\tau_0<t-s)$.
Thus if
$\PP^r_0(\tau^r_0<t-s,\hat X^r(\tau^r_0)\in B^r_i)\to q_i\PP_0(\tau_0<t-s)$
then we obtain
\[
\PP^r_{y^r}(\tau^r_m\le s,\chi^r_m\in B^r_i)-q_i\PP^r_{y^r}(\tau^r_m\le s)
\to0,
\]
and since $\PP^r_{y^r}(\tau^r_m\le s)\to \PPrbm_{\|x\|}(\tau_m\le s)$ uniformly over $u$-admissible pairs $(x,y)$, \eqref{31} follows.

To prove \eqref{31} for $y\equiv0$ and $m=0$, note that
under $\PP^r_0$, $\zeta^r_0=0$ and so $\tau^r_0$ is a.s.\ equal to
$\tau^r=\tau^r(\eps)=\inf\{s\ge 0:\hat R^r(s)\ge \eps\}$ (see \eqref{35}).
Moreover, $\chi^r_0$
that has been defined as $\hat X^r(\tau^r_0)$ is a.s.\ equal to $\chi^r:=\hat X^r(\tau^r)$.
Hence we aim now at showing
\begin{equation}
  \label{32}
  \PP^r_0(\tau^r\le s,\chi^r\in B^r_i)\to \PPrbm_0(\tau\le s)q_i.
\end{equation}
Without loss of generality, we take $i=1$.
In addition to the parameter $\eps$, that has been fixed,
we introduce a new parameter,
$a\in(0,\eps)$, that will play the role of the parameter $\eps$ in Proposition \ref{lem3}.
We introduce several pieces of notation associated with $a$ in a way
analogous to those defined in terms of $\eps$. Namely, $B^{r,a}_i=B(ae_i,u(r))$,
$\tau^{r,a}=\inf\{s\ge 0:\hat R^r(s)\ge a\}$ and $\chi^{r,a}=\hat X^r(\tau^{r,a})$.
In addition, we let $\nu^{r,a}_i$ denote the probability measures
supported on $B^{r,a}_i$, given by $\PP^r_0(\chi^{r,a}\in\cdot|\chi^{r,a}\in B^{r,a}_i)$.

Let
\[
g^r(z)=\PP^r_{z}(\tau^r<t,\chi^r\in B^r_1).
\]
We analyze $g^r(0)$ by studying its relation to
$g^{r,a}_i:=\int g^r(z)\nu^{r,a}_i(dz)$.
First,
\begin{align}\label{33}
  g^r(0)&=\E^r_0[\E^r_0[\one_{\{\tau^r<t,\chi^r\in B^r_1\}}|\calF^r_{\tau^{r,a}}]]
  = \E^r_0[\psi^r(\tau^{r,a},\chi^{r,a})],
\end{align}
where
\[
\psi^r(s,z)=\PP^r_{z}(\tau^r<t-s,\chi^r\in B^r_1).
\]
Hence
\begin{align}\label{33a}
g^r(0)=\E^r_0[\psi^r(0,\chi^{r,a})]+\del^r_a,
\end{align}
where
\[
|\del^r_a|\le \E^r_0[|\psi^r(\tau^{r,a},\chi^{r,a})-\psi^r(0,\chi^{r,a})|].
\]

We denote by $O^r_a(h(a))$ (resp., $o^r_a(h(a))$)
any function $g$ of the tuple $(x,y,r,a)$ satisfying
$\limsup_{a\down0}\limsup_{r\to\iy}\sup_{(x,y)\text{ $u$-admissible}}h(a)^{-1}|g(x,y,r,a)|<\iy$
(resp., $=0$).
We argue that $\del^r_a=O^r_a(a^2)$. To this end, note that
\begin{align*}
  0\le\psi^r(0,z)-\psi^r(s,z)
  &=\PP^r_z(\tau^r<t,\chi^r\in B^r_1)-\PP^r_z(\tau^r<t-s,\chi^r\in B^r_1)\\
  &\le \PP^r_z(t-s<\tau^r<t).
\end{align*}
Now, $\tau^r\To\tau$ as $t\to\iy$, uniformly for $z$ in $B(0,\eps/2)$.
Moreover, for RBM, denoting the density $\frac{d}{d\theta}\PPrbm_\eta(\tau\le\theta)$
by $l(\eta,\theta)$ (where, as before, $\tau=\tau(\eps)$), a uniform bound
holds in the form
\begin{equation}\label{110}
\sup_{\eta\in[0,\eps/2],\theta\in[t/2,\iy)}l(\eta,\theta)<\iy.
\end{equation}
Indeed, an explicit eigenfunction expansion of the density $l$ is given in \cite{hitting}.
Using equations (3.15)--(3.19) of \cite{hitting} one can directly
obtain the bound $\sup_{x\in[0,\eps],t\ge t_0}l(x,t)<\iy$
for any constant $t_0>0$. This gives \eqref{110}.
Using \eqref{110} for $s\in[0,t/2]$ and the trivial bound $1$ for
$s\in[t/2,t]$ gives
\[
\sup_{\eta\in[0,\eps/2]}\PPrbm_\eta(t-s<\tau<t)\le cs,
\]
for some constant $c$ (which may depend on $t$), for all $s\in[0,t]$.
In view of this, $\limsup_r|\del^r_a|\le c\limsup_r\E^r_0[\tau^{r,a}]\le ca^2$,
where the last inequality is standard, and follows by Brownian scaling.

Going back to \eqref{33a} and noting that $\psi^r(0,z)=g^r(z)$, we have
$g^r(0)=\E^r_0[g^r(\chi^{r,a})]+O^r_a(a^2)$.
Therefore, it follows from Lemma \ref{lem-b}(ii) that
the probability of having $\chi^{r,a}\notin \cup_iB^{r,a}_i$ is $o_r(1)$, hence
\begin{align}
\notag
  \label{26}
  g^r(0)&=\sum_i \PP^r_0(\chi^{r,a}\in B^{r,a}_i)g^{r,a}_i+O^r_a(a^2)\\
  &=\sum_i q^{r,a}_ig^{r,a}_i+O^r_a(a^2),
\end{align}
where by Proposition \ref{lem3}, $q^{r,a}_i:=\PP^r_0(\chi^{r,a}\in B^{r,a}_i)=q_i+o^r_a(1)$
(note that $B^{r,a}_i$ agrees with the ball from Proposition \ref{lem3}).

Next consider initial condition $\nu^{r,a}_i$, for which we can write
\begin{align*}
g^{r,a}_i
&=\E^r_{\nu^{r,a}_i}[g^r(\hat X^r(0))]
=\E^r_{\nu^{r,a}_i}[\hat\psi^r(\tilde\tau^r,\hat X^r(\tilde\tau^r))],
\end{align*}
where $\tilde\tau^r=\inf\{s\ge 0:\hat R^r(s)\notin(0,\eps)\}$ and $\hat\psi^r(s,z)=\PP^r_z(\tau^r<t-s,\chi^r\in B^r_1)$.
Now,
\begin{align*}
\hat\psi^r(s,0)&=\PP^r_0(\tau^r<t-s,\chi^r\in B^r_1),
\end{align*}
and for every $s$ that satisfies $\rho(s)\ge \eps$,
\begin{align*}
\hat\psi^r(s,z)=\one_{\{s\le t\}},\;\text{for } z\in B^r_1\quad\text{and}\quad
\hat\psi^r(s,z)=0,\;\text{for } z\in B^r_i,\ i\ne1.
\end{align*}
Hence
\begin{align}\label{27}
\notag
g^{r,a}_i&=\E^r_{\nu^{r,a}_i}[\one_{\{\hat R^r(\tilde\tau^r)=0\}}\hat\psi^r(\tilde\tau^r,0)]
+\E^r_{\nu^{r,a}_i}[\one_{\{\hat R^r(\tilde\tau^r)\ge \eps\}}\hat\psi^r(\tilde\tau^r,\hat X^r(\tilde\tau^r))]\\
&=\PP^r_{\nu^{r,a}_i}(\hat R^r(\tilde\tau^r)=0)\hat\psi^r(0,0)
+\one_{\{i=1\}}\PP^r_{\nu^{r,a}_i}(\hat R^r(\tilde\tau^r)\ge \eps,\tilde\tau^r<t)+\hat\del^r_a+o_r(1),
\end{align}
where
\[
|\hat\del^r_a|\le \E^r_{\nu^{r,a}_i}[\one_{\{\hat R^r(\tilde\tau^r)=0\}}|\hat\psi^r(\tilde\tau^r,0)-\hat\psi^r(0,0)|],
\]
and we used again the bound \eqref{13}.
To further bound $\hat\del^r_a$, note that the argument provided earlier
for $\psi^r$ can be used also for $\hat\psi^r$, and gives
$|\hat\del^r_a|\le c\E^r_{\nu^{r,a}_i}[\one_{\{\hat R^r(\tilde\tau^r)=0\}}\tilde\tau^r]$.
On the indicated event, $\tilde\tau^r$ is bounded by the exit time of
$\hat R^r$ from the interval $(0,2a)$,
the expectation of which is $O^r_a(a^2)$. Hence $\hat\del^r_a=O^r_a(a^2)$.

Next, for the RBM $\rho$ denote analogously $\tilde\tau=\inf\{s\ge0:\rho(s)\notin(0,\eps)\}$.
Denote
\[
\beta_1(a)=1-\PPrbm_{a}(\rho(\tilde\tau)=0)=\PPrbm_{a}(\rho(\tilde\tau)=\eps),
\qquad
\beta_2(a)=\PPrbm_a(\rho(\tilde\tau)=\eps,\tilde\tau\le t).
\]
Then $\PP^r_{\nu^{r,a}_i}(\hat R^r(\tilde\tau^r)=0)=1-\beta_1(a)+o_r(1)$
for all $i$. Moreover, $\PP^r_{\nu^{r,a}_i}(\hat R^r(\tilde\tau^r)\ge \eps,\tilde\tau^r<t)
=\beta_2(a)+o_r(1)$.
Hence from \eqref{26} and \eqref{27} we obtain
\[
\begin{cases}
  g^r(0)=\ds\sum_i(q_i+o^r_a(1))g^{r,a}_i+O^r_a(a^2),
  \\
  g^{r,a}_1=(1-\beta_1(a)+o_r(1))g^r(0)+\beta_2(a)+O^r_a(a^2)+o_r(1),
  \\
  g^{r,a}_i=(1-\beta_1(a)+o_r(1))g^r(0)+O^r_a(a^2)+o_r(1),\qquad\qquad i\ne1.
\end{cases}
\]
Solving this system of equations gives
\[
g^r(0)=\frac{(q_1+o^r_a(1))\beta_2(a)+O^r_a(a^2)}{\beta_1(a)+o_r(1)}.
\]
Denote
\[
G_a=\frac{(q_1+o_a(1))\beta_2(a)+O_a(a^2)}{\beta_1(a)}.
\]
In order to show that $\lim_rg^r(0)=q_1\PP_0(\tau<t)$ it suffices to show
that $\lim_{a\downarrow0}G_a=q_1\PP_0(\tau<t)$.
Since it is known for RBM (equivalently, for a 1-dimensional BM)
that $\beta_1(a)>ca$ for some constant $c>0$
and all small $a$, it suffices to show that $\beta_2(a)/\beta_1(a)\to \PPrbm_0(\tau<t)$
as $a\down0$. To this end, use strong Markovity to write
\[
\PPrbm_a(\rho(\tilde\tau)=0,\tau\le t)
=\Erbm_a[\one_{\{\rho(\tilde\tau)=0\}}\ph^{\#}(\tilde\tau,\rho(\tilde\tau))],
\]
$\ph^{\#}(s,x)=\PPrbm_x(\tau\le t-s)$ for $x\in\R_+$. Now, $0\le \ph^\#(0,0)-\ph^\#(s,0)\le cs$,
and therefore
\[
|\PPrbm_a(\tau\le t,\rho(\tilde\tau)=0)-\PPrbm_0(\tau\le t)\PPrbm_a(\rho(\tilde\tau)=0)|
\le c\E_a[\one_{\{\rho(\tilde\tau)=0\}}\tilde\tau]\le ca^2.
\]
Hence
\begin{align*}
\beta_2(a)
&=\PPrbm_a(\tau\le t)-\PPrbm_a(\tau\le t,\rho(\tilde\tau)=0)
\\
&=\PPrbm_a(\tau\le t)-\PPrbm_0(\tau\le t)\PPrbm_a(\rho(\tilde\tau)=0)+O(a^2)
\\
&=\PPrbm_0(\tau\le t)\beta_1(a)+\del^\#(a)+O(a^2),
\end{align*}
where $\del^\#(a)=\PPrbm_a(\tau\le t)-\PP_0(\tau\le t)$.
If we show that $\del^\#(a)=O(a^2)$ then
$\beta_2(a)/\beta_1(a)\to \PPrbm_0(\tau\le t)$ as $a\down0$, and the proof is established.

To show that $\del^\#(a)=O(a^2)$,
let $\calL$ denote the generator of the process $\rho$, and
$v(x,s)=\PPrbm_x(\tau>s)$ for $x\in\R_+$. Then $\calL$ is given by
$\frac{\sig^2}{2}\frac{d^2}{dx^2}+b\frac d{dx}$,
with the Neumann boundary condition at $0$ and the Dirichlet boundary condition
at $\eps$, and $v$ is a smooth function satisfying
\begin{align*}
  &\partial_sv=\calL v, \hspace{6.5em} x\in(0,\eps),\ s>0,\\
  &\partial_xv(0,s)=v(\eps,s)=0, \hspace{1em} s>0,\\
  &v(x,0)=1, \hspace{5.2em} x\in(0,\eps).
\end{align*}
In particular, it is a smooth function satisfying $\partial_xv(0,s)=0$,
and therefore $v(x,t)=v(0,t)+O(x^2)$, for $t$ fixed
and $x\down0$. This shows that $\del^\#(a)=O(a^2)$.
\qed

\subsection{The total nominal workload process}\label{sec32}

In this section we prove Lemma \ref{lem1}.
Roughly stated, this lemma asserts that the total nominal workload process
converges at diffusion scale to an RBM. This is a well-understood
fact for an {\it arbitrary} non-idling policy. However, for completeness
and since the statement of the lemma involves
uniform convergence, which is perhaps less standard, we provide a proof.

\skp
\noi
\noi
{\bf Proof of Lemma \ref{lem1}.}
i. We start the proof with some notation aimed at describing the scaled nominal workload process
in terms of scaled arrival and service processes.
Let $\bar T_i(t)=\frac{\la_i}{\mu_i}t$, and
\begin{align}\notag
	\hat A^r_i(t)&=r^{-1}(A^r_i(t)-\la^r_i t),
	\hspace{3.5em}
	\hat S^r_i(t)=r^{-1}(S^r_i(t)-\mu^r_i t)
	,\\\notag
	\hat Y^r_i(t)&=\mu^r_ir^{-1}\Big(\bar T_i(t)-T^r_i(t)\Big),\qquad b^r_i= r^{-1}\Big(\la^r_i-\frac{\la_i}{\mu_i}\mu^r_i\Big),
\end{align}
for $t\in\R_+$. Then by \eqref{38},
\begin{align}\label{as01}
	\hat Q^r_i(t)=\hat Q^r_i(0)+\hat A^r_i(t)-\hat S^r_i( T^r_i(t))+b^r_i t+\hat Y^r_i(t),\qquad t\in\R_+.
\end{align}
Note by \eqref{asaf1} and \eqref{14} that $\hat R^r=\one\cdot\hat X^r
=\sum_i\frac{r^2}{\mu^r_i}\hat Q^r_i$.
If we denote
\begin{align*}
B^r(t)&=\sum_i\frac{r^2}{\mu_i^r}[\hat A^r_i(t)-\hat S^r_i( T^r_i(t))+b^r_it],
\\
U^r(t)&=\sum_i\frac{r^2}{\mu_i^r}\hat Y^r_i(t)=r[t-\one\cdot T^r(t)],
\end{align*}
then we have the identity
$\hat R^r=B^r+U^r$. Moreover, by its definition, $\hat R^r$
takes values in $\R_+$. Furthermore, by the non-idling property,
the right derivative of $\one\cdot T^r$ at $t$ assumes the value 1
if and only if the system is non-empty at this time,
that is, $\hat R^r(t)>0$ (it otherwise assumes the value 0).
Consequently, $\int_0^\iy\hat R^r(t)dU^r(t)=0$. These three properties
imply
\begin{align}\label{eq_RBM}
	(\hat R^r,U^r)=\Gam[\hat R^r(0)+B^r].
\end{align}
It follows from the expression \eqref{40} for $\Gamma$ that
for $t>0$
\begin{align}\label{1000}
|\hat R^r(t)-\hat R^r(0)|\le 2\|B^r\|_t.
\end{align}
Now, the bound $T^r_i(t)\le t$ for all $t$ gives
$\|\hat S^r_i\circ T^r_i\|_T\le\|\hat S^r_i\|_T$.
The quantities $r^2/\mu^r_i$ as well as $b^r_i$ converge
in view of \eqref{204}. Hence $\|B^r\|_T\le c(\|\hat A^r\|_T+\|\hat S^r\|_T+1)$.
Therefore,
\[
\|\hat R^r-\hat R^r(0)\|_T\le
	c(\|\hat A^r\|_T+\|\hat S^r\|_T+1).
\]
By the functional central limit theorem for renewal processes (see Theorem 14.6 of \cite{Bill}),
$(\hat A^r,\hat S^r)$ converge to a BM with drift zero and diffusion matrix $\Xi$,
where
\[
\hat A^r=(\hat A^r_i)_{i=1}^N,\quad \hat S^r=(\hat S^r_i)_{i=1}^N,\quad \Xi=\text{diag}(\la_1^{1/2},\ldots,\la_N^{1/2},\mu_1^{1/2},\ldots,\mu_N^{1/2}).
\]
This implies that $\|\hat A^r\|_T+\|\hat S^r\|_T$ is a tight sequence of RVs
(for each fixed $T$), and in view of the above bound, so is $\|\hat R^r\|_T$.

By the discussion preceding Theorem \ref{th1}, $\hat X^r$ and $\hat Q^r$
are asymptotically related via the matrix $\hat M$. Appealing to \eqref{as01} again and recalling that $\one\cdot \hat X^r=\hat R^r$, we have $\|\hat Y^r\|_T\le c(\|\hat A^r\|_T+\|\hat S^r\|_T+1)$, and we get the tightness of $\|\hat Y^r\|_T$ (uniformly in the initial state).
In view of the definition of $\hat Y^r$, we obtain that for every $T,v>0$,
\begin{align}\label{cohen3}
\limsup_{r\to\iy}\sup_{x\in\calS^r}\bP^r_x(\|T^r-\bar T\|_T>v)=0.
\end{align}
Set
\begin{align*}
\tilde B^r(t)&=\sum_i\frac{r^2}{\mu_i^r}[\hat A^r_i(t)-\hat S^r_i( \bar T_i(t))+b_it],\qquad E^r=B^r-\tilde B^r.
\end{align*}
Notice that $\tilde B^r$ is measurable w.r.t.~$\sigma\{A^r(t),S^r(t), t\in\R_+\}$.
Now,
\[
\|E^r\|_T\le\sum_i\frac{r^2}{\mu^r_i}\|\hat S^r_i\circ T^r_i-\hat S^r_i\circ\bar T_i\|_T
\le c\sum_iw_T(\hat S^r_i,\theta^r),
\]
where $\theta^t=\|T^r-\bar T\|_T$. The $\calC$-tightness of $\hat S^r$
along with \eqref{cohen3} give \eqref{cohen1}.
Finally, the convergence in law of $(\hat A^r,\hat S^r)$
gives $\tilde B^r\To B$, where $B$ is a $(b,\sig)$-BM.

\skp
\noi ii. Fix $v,t_0>0$.
For every $x\in\calS^r$ with $\one\cdot x<v$,
one has by the representation $\hat R^r=\Gam_1[\hat R^r(0)+B^r]$,
noting that $B^r$ is measurable w.r.t.\ $\sig\{A^r,S^r,X^r\}$,
\begin{align}\notag
	\PP^r_x(\zeta^r\le t_0)&=\PP^r_x\Big(\inf_{t\in[0,t_0]}\hat R^r(t)=0\Big)\\\notag
	&=\bP^r_x\Big(\hat R^r(0)+\inf_{t\in[0,t_0]}B^r(t)\le 0\Big)\\\notag
	&\ge \bP^r_x\Big(\inf_{t\in[0,t_0]}B^r(t)<-v\Big).
\end{align}
By part (i) of the lemma, specifically, the convergence of $\tilde B^r$ to $B$ (indep of $x$)
and the uniform estimate \eqref{cohen1} on $E^r$,
it follows that for every $\delta>0$ and all sufficiently large $r$
and $x\in\calS^r$ with $\one\cdot x\le v$, the RHS of the above display is bounded below
by $\bP(\inf_{t\in[0,t_0]}B(t)<-2v)-\delta$. The last expression does not depend on $x$
and, as $B$ is a BM starting at the origin, converges to $1-\del$ as $v\down0$.
Therefore
\begin{align}\notag
	\liminf_{\del\down0}\liminf_{r\to\iy}\inf_{x\in\calS^r:\one\cdot x<v}\PP^r_x(\zeta^r\le t_0)
	\ge1-\del.
\end{align}
Taking $\del\down0$ gives the result.
\qed

\subsection{Estimates on exiting the tubes}\label{sec33}

In this section we develop an estimate on the displacement of the prelimit
process $\hat X^r$ away from $\calS_0$.
The main use of this estimate is in the argument provided in \S\ref{sec34}.
In addition, the statement constitutes a strong form of that of Lemma \ref{lem-b}(i).
Thus at the end of the section we provide a proof of Lemma \ref{lem-b} based
on this estimate.

The proof is based on a Lyapunov function technique.
This function is constructed so that it expresses
the total nominal workload in all buffers save the one where queue length is greatest.
For a precise definition we need some notation. Recall from \eqref{39}
the sets $\calK(x)$, and note that
$\calK(\hat X^r(t))$ gives the set of shortest nonempty queues at time $t$.
For $x\in\R_+^N$, let
\[
\calM(x)=\{i\in[N]:x_i\ge x_j \text{ for all } j\in[N]\}.
\]
Then $\calM(\hat X^r)$ gives the set of longest queues.
Let $F:\R_+^N\to\R$ be given by
\begin{align}\label{eq100}
F(x)=\sum_ix_i-\max_ix_i.
\end{align}
Note that $F$ is nonnegative and vanishes on the set $\calS_0$ and only there.
\begin{lemma}\label{lem2}
Given $c_0>0$, $\kap_0\in(0,1/2)$ and $0<\gamma_1<\gamma_2<\iy$,
there exist constants $r_0, c_1>0$ such that for every $r>r_0$ and every initial state
$x\in\calS^r$ that satisfies $F(x)\le \gamma_1 r^{-\kap_0}$,
\begin{align}\label{eq99}
\PP^r_x\left(\|F(\hat X^r(\cdot))\|_{c_0\log r}>\gamma_2r^{-\kap_0}\right)\le r^{-c_1}.
\end{align}
\end{lemma}

Lemma \ref{lem2} and the first item of Lemma \ref{lem-b} are similar, where the former
is concerned with long time intervals
as well as rates of convergence. However, the latter is not an immediate
consequence of the former. We present their proofs together.

\skp

\noi {\bf Proof of Lemma \ref{lem2} and Lemma \ref{lem-b}(i).}
For the proof of Lemma \ref{lem2},
fix $c_0$, $\kap_0$, $\gamma_1$ and $\gamma_2$ as in the statement of the lemma.
Using the expression \eqref{10} for the generator of $X^r$
write the one for $\hat X^r=rX^r$ (see \eqref{14}), as
\begin{equation}
  \label{10+}
  \calL^rf(x)=\sum_i\la^r_i\Big(f\Big(x+\frac{r}{\mu^r_i}e_i\Big)-f(x)\Big)
  +\sum_{i\in\calK(x)}p^{\calK(x)}_i\mu^r_i\Big(f\Big(x-\frac{r}{\mu^r_i}e_i\Big)-f(x)\Big),
\end{equation}
for bounded $f:\calS^r\to\R$.

Recall that $c$ denotes a generic positive constant that
does not depend on $r$. We begin by showing that there exists a constant $c$
such that for all $r$ sufficiently large,
\begin{align}\label{eq103k}
{\cal L}^rF(x)< -c r \text{ for all $x$ such that $F(x)>0$}.
\end{align}
To this end, note that
the first term on the RHS of \eqref{10+}, upon substituting $F$ for $f$, equals
\begin{align*}
&\sum_i \la^r_i\Big(\frac{r}{\mu^r_i}+\max\{x_1,\ldots,x_N\}
-\max\{x_1,\ldots,x_{i-1},x_i+r/\mu^r_i,x_{i+1},\ldots,x_N\}\Big)
\\
&\le r\sum_{i\in[N]\setminus\calM(x)}\frac{\la^r_i}{\mu^r_i}.
\end{align*}
The inequality above is valid since for $i\in\calM(x)$ the $i$-th term in the sum is zero,
and for $i\notin \calM(x)$,
\[
\max\{x_1,\ldots,x_N\}
\le\max\{x_1,\ldots,x_{i-1},x_i+r/\mu_i^r,x_{i+1}\ldots,x_N\}.
\]
The second term on the RHS of \eqref{10+} (with $f=F$) can be expressed as
\begin{equation}\label{36}
\sum_{i\in \calK(x)}p^{\calK(x)}_i\mu^r_{i}
\Big(-\frac{r}{\mu^r_i}+\max\{x_1,\ldots,x_N\}
-\max\{x_1,\ldots,x_{i-1},x_i-r/\mu_i^r,x_{i+1},\ldots,x_N\}\Big).
\end{equation}
We argue that for $i\in \calK(x)$,
\[
\max\{x_1,\ldots,x_N\}
=\max\{x_1,\ldots,x_{i-1},x_i-r/\mu_i^r,x_{i+1}\ldots,x_N\}.
\]
If $\calK(x)=\calM(x)=\{x_j\}$ for some $j\in[N]$ then $F(x)=0$. Therefore, if $F(x)>0$
then either $\calK(x)\ne\calM(x)$ or $\calK(x)=\calM(x)$ and $\calK(x)$ contains more than one element.
In both cases, for every $i\in \calK(x)$ there is $j\in\calM(x)$,
different from $i$, such that both maxima above equal $x_j$.
This shows that the expression in \eqref{36} equals $-r$.
Combining this with the bound on the first term,
and recalling that $\la^r_i/\mu^r_i$ is asymptotic to $\la_i/\mu_i$,
and that the latter fractions sum to 1, shows \eqref{eq103k}.

We analyze the event $\Om^r:=\{\|F(\hat X^r(\cdot))\|_{c_0\log r}>\gamma_2r^{-\kap_0}\}$
under $\PP^r_x$ for $x$ such that $F(x)\le \gamma_1 r^{-\kap_0}$.
Recall that the jump sizes of $\hat X^r$ are at the scale of $r^{-1}$; as a result,
the same is true for the process $F(\hat X^r(\cdot))$.
Since $\kap_0<1/2$, $r^{-\kap_0}$ is at a larger scale than these jumps. Hence
there exist random times $0\le \theta^r_1<\theta^r_2\le c_0\log r$ such that $\PP^r_x$-a.s.\ on $\Om^r$,
\begin{align}\label{eq101}
F(\hat X^r(\theta^r_1))\le  \gamma_1 r^{-\kap_0},
\quad F(\hat X^r(\theta^r_2))\ge \gamma_2r^{-\kap_0},
\quad \text{and}\quad 0<F(\hat X^r(t))\le\gamma_2r^{-\kap_0},
\quad t\in[\theta^r_1,\theta^r_2).
\end{align}
The process
\begin{align}\label{2000}
M^r(t)=F(\hat X^r(t))-F(x)-\int_0^t{\cal L}^rF(\hat X^r(s))ds,\qquad t\in\R_+,
\end{align}
is a local martingale.
From \eqref{eq103k} and \eqref{eq101}, denoting $\del=\gamma_2-\gamma_1>0$,
one has
\begin{align}\label{eq105}
M^r(\theta^r_2)-M^r(\theta^r_1)\ge c r(\theta^r_2-\theta^r_1)+\del r^{-\kap_0}.
\end{align}
Fix a constant $d\in(2\kap_0,1)$ and consider the events
\[
\Om^r_1=\{\text{$\theta^r_2-\theta^r_1\le r^{-d}$ and \eqref{eq105} holds}\},
\qquad
\Om^r_2=\{\text{$\theta^r_2-\theta^r_1>r^{-d}$ and \eqref{eq105} holds}\}.
\]
Then $\PP^r_x(\Om^r)\le \PP^r_x(\Om^r_1)+\PP^r_x(\Om^r_2)$.
We argue separately for the two events.

\skp
\noi{\it The event $\Om^r_1$.}
Let intervals $I^r_j$ be defined by $I^r_j=[jr^{-d}/2,(j+1)r^{-d}/2]$ for
$j\in\{0,1,\ldots,j_1(r)\}$, where $j_1(r)=\lfloor 2c_0r^d\log r\rfloor$.
On $\Om^r_1$ there must exist $j$
and an interval $J\subset I^r_j\cap(\theta^r_1,\theta^r_2)$ such that
\[
\osc_{J}M^r \ge \del r^{-\kap_0}/3,
\]
where here and throughout, $\osc_Af=\sup_Af-\inf_Af$.
As a result, $\osc_{I^r_j}M^r\ge\del r^{-\kap_0}/3$.
Therefore, using the Burkholder-Davis-Gundy (BDG) inequality \cite[Theorem 48]{Protter2004}
and denoting by $[ M^r]_I$ the quadratic variation of $M^r$ over an interval $I$,
\begin{align}\label{eq102}
\PP^r_x(\Om^r_1)
\le\sum_{j=0}^{j_1(r)}
\PP^r_x\Big(\osc_{I^r_j}M^r\ge \frac{\del r^{-\kap_0}}{3}\Big)
\le C_k\sum_{j=0}^{j_1(r)}
\Big(\frac{6r^{\kap_0}}{\del}\Big)^{2k}\E^r_x\{[M^r]^k_{I^r_j}\},
\end{align}
where $k$ is any number in $[1/2,\iy)$.
The quadratic variation process $[M^r]$ has piecewise-constant samples paths
with jumps taking values in the set $\{(r/\mu^r_1)^2,\ldots,(r/\mu^r_N)^2\}$.
The number of its jumps in an interval $[s,t]$ is stochastically dominated
by a Poisson RV with parameter $(t-s)\sum_i(\mu^r_i+\la^r_i)$.
Since $\mu^r_i$ and $\la^r_i$ scale like $r^2$, $[M^r]_{I^r_j}$
is stochastically dominated by $K^r=cr^{-2}\pi^r$, where
$\pi^r$ is a Poisson RV with parameter $cr^{2-d}$.
Consequently,
\[
\E^r_x\{[M^ir]^k_{I^r_j}\}\le cr^{-dk}.
\]
Therefore the r.h.s.~of \eqref{eq102} is bounded above by
$cr^{d-k(d-2\kap_0)}\log r$ (where $c$ may depend on $k$).
Taking $k>\max\{1/2,d/(d-2\kap_0)\}$
gives the bound $\PP^r_x(\Om^r_1)\le r^{-c}$, provided that $r$ is sufficiently large.

\skp
\noi{\it The event $\Om^r_2$.}
Clearly,
\[
\PP^r_x(\Om^r_2)\le \PP^r_x(\osc_{[0,c_0\log r]}M^r\ge cr^{1-d}).
\]
Using again the BDG inequality (with $k=1/2$)
followed by a domination of the number of jumps in terms of a Poisson RV
with parameter $c_0r^2\log r$, and the sizes of the jumps by $r^{-2}$, gives
\begin{align}\notag
\PP^r_x(\Om^r_2)\le\PP^r_x(2\|M^r\|_{c_0\log r}\ge cr^{1-d})
\le c\frac{\E^r_x\left\{[ M^r]_{c_0\log r}\right\}}{r^{2(1-d)}}
\le \frac{c\log r}{r^{2(1-d)}}<r^{-c}.
\end{align}
This completes the proof of Lemma \ref{lem2}.

\skp
In order to establish the proof of Lemma \ref{lem-b}(i) we prove below the following stronger result that also serves us in the proof of Lemma \ref{lem-b}(iii): for every $u\in U_0$ satisfying $\lim_{r\to\iy}r(u(r))^3=\iy$, one has
\begin{align}\label{COHEN002}
\liminf_{r\to\iy}
\inf_{x\in\calS^r\cap\calS_{u(r)/2}}
\PP^r_x(\hat X^r(t)\in\calS_{u(r)} \text{ for all } t\in[0,T])=1.
\end{align}
This statement implies Lemma \ref{lem-b}(i), since by the assumption on $\nu^r$, there exists $ u\in\calU_0$,
such that $\nu^r(\hat X^r(0)\in\calS_{ u(r)/2})\to1$;
without loss of generality we may assume that $\lim_{r\to\iy}r(u(r))^3=\iy$.

We next show how the details of the above proof are modified in order to
prove \eqref{COHEN002}. Fix an arbitrary $ u\in\calU_0$ that satisfies
$\lim_{r\to\iy}r(u(r))^3=\iy$. We claim that
\begin{align}\notag
\limsup_{r\to\iy}
\sup_{x\in\calS^r\cap\calS_{u(r)/2}}
\PP^r_{x}\left(\|F(\hat X^r(\cdot))\|_{T}>u(r)\right)=0.
\end{align}
Unlike in \eqref{eq99}, we consider here a fixed horizon $T$, we do not provide a convergence rate, and the polynomial tube widths are replaced by $u(r)$.

Recall \eqref{10+}, \eqref{eq103k}, and \eqref{2000}.
We analyze the event $\bar\Om^r:=\{\|F(\hat X^r(\cdot))\|_{T}>u(r)\}$
under $\PP^r_{x}$.
Since $ru(r)\to\iy$ and the jump sizes of the process $F(\hat X^r(\cdot))$ are at the scale of $r^{-1}$, there must exist random times $0\le \bar\theta^r_1<\bar\theta^r_2\le T$ such that $\PP^r_{x}$-a.s.\ on $\bar\Om^r$,
\begin{align}\label{eq101new}
F(\hat X^r(\bar\theta^r_1))\le  \frac{1}{2}u(r),
\quad F(\hat X^r(\bar\theta^r_2))\ge u(r),
\quad \text{and}\quad 0<F(\hat X^r(t))\le u(r),
\quad t\in[\bar\theta^r_1,\bar\theta^r_2).
\end{align}
From \eqref{eq103k} and \eqref{eq101new},
one has
\begin{align}\label{eq105new}
M^r(\bar\theta^r_2)-M^r(\bar\theta^r_1)\ge c r(\bar\theta^r_2-\bar\theta^r_1)+\frac{1}{2}u(r).
\end{align}
Set $u_0(r)=(u(r))^3$ and consider the events
\[
\bar\Om^r_1=\{\text{$\bar\theta^r_2-\bar\theta^r_1\le u_0(r)$ and \eqref{eq105new} holds}\},
\qquad
\bar\Om^r_2=\{\text{$\bar\theta^r_2-\bar\theta^r_1>u_0(r)$ and \eqref{eq105new} holds}\}.
\]
Then $\PP^r_{x}(\bar\Om^r)\le \PP^r_{x}(\bar\Om^r_1)+\PP^r_{x}(\bar\Om^r_2)$.
We argue separately for the two events.

\skp
\noi{\it The event $\bar\Om^r_1$.}
Let intervals $\bar I^r_j$ be defined by
\[
\bar I^r_j=[ju_0(r)/2,(j+1)u_0(r)/2]
\]
for
$j\in\{0,1,\ldots,\lfloor 2T/u_0(r)\rfloor\}$. The same arguments given before with the choice of $k=2$ in BDG inequality lead to the following sequence of inequalities and the uniform limit over $\calS^r\cap\calS_{u(r)/2}$
\begin{align}\label{eq102new}
\begin{split}\PP^r_{x}(\Om^r_1)
&\le\sum_{j=0}^{j_1(r)}
\PP^r_{x}\Big(\osc_{I^r_j}M^r\ge \frac{u(r)}{6}\Big)
\le C_k\sum_{j=0}^{j_1(r)}
\Big(\frac{12}{u(r)}\Big)^{4}\E^r_{x}\{[M^r]^2_{I^r_j}\}\\
&\le
\sum_{j=0}^{j_1(r)}
\Big(\frac{12}{u(r)}\Big)^{4}(u_0(r))^2\le cu^2(r)\to0.
\end{split}
\end{align}

\skp
\noi{\it The event $\bar \Om^r_2$.} Arguing as before, we obtain,
\begin{align}\notag
\PP^r_{x}(\Om^r_2)\le\PP^r_{x}(2\|M^r\|_T\ge cru_0(r))
\le c\frac{\E^r_{x}\left\{[ M^r]_{T}\right\}}{(ru_0(r))^2}
\le \frac{c}{(ru_0(r))^2}.
\end{align}
By our choice of the function $u$, the last expression converges to $0$ as $r\to\iy$, uniformly
for $x\in\calS^r\cap\calS_{u(r)/2}$.
\hfill$\Box$

\skp

\noi{\bf Proof of Lemma \ref{lem-b} (continued).}
First, the assertion regarding $\calC$-tightness follows directly from Lemma \ref{lem1}(i) and Lemma \ref{lem-b}(i).

\skp
\noi ii. The statement of this part follows from part (i) with initial condition $0$ and strong Markovity.

\skp\noi
iii.
It is sufficient to show that for every $u\in\calU_0$ satisfying $\lim_{r\to\iy}r(u(r))^3=\iy$, one has
\begin{equation}\label{90}
\limsup_{r\to\iy}
\sup_{y\in\calS^r,x\in[0,k]:\|y-xe_i\|<u(r)/2}
\Big|\E^r_y[f(\hat X^r(t))\one_{\{t<\zeta^r\}}]
-\Erbm_x[f(\rho(t)e_i)\one_{\{t<\zeta\}}]\Big|=0,
\end{equation}
We first show that for every such $u$ and every $\eps>0$,
\[
\limsup_{r\to\iy}
\sup_{y\in\calS^r,x\in[0,k]:\|y-xe_i\|<u(r)/2}
\PP^r_y\left(\zeta^r>t,\hat R^r(t)>\eps,\|\hat X^r(t)-\hat X^r_i(t)e_i\|>u(r)\right)=0.
\]
Indeed, in order for the process $\hat X^r$, starting inside the $u(r)/2$-tube around axis $i$,
to exit the $u(r)$-tube around the same axis by time $t$ and reach $c\eps$ away from the origin,
it must either escape $\calS_{u(r)}$ before $t$
or pass through a $cu(r)$-neighborhood of the origin without hitting
the origin and then move through a different tube and $c\eps$ away from the origin.
The probabilities of these two events converge to zero uniformly
in the initial conditions; the first convergence follows by
\eqref{COHEN002}. The second event can be expressed in terms of
an atypical behavior of $\hat R^r$,
as a sequence of processes converging in law to an RBM.

The assertion above along with the uniform continuity of $f$
and a further application of \eqref{COHEN002}, imply
\[
\limsup_{r\to\iy}
\sup_{y\in\calS^r,x\in[0,k]:\|y-xe_i\|<u(r)/2}
\Big|\E^r_y[f(\hat X^r(t))\one_{\{t<\zeta^r\}}]
-\E^r_y[f(\hat R^r(t)e_i)\one_{\{t<\zeta^r\}}]\Big|=0.
\]
Finally, the statement in Lemma \ref{lem1}(i),
according to which $\hat R^r\To\rho$ holds with an error term that converges to zero
uniformly in the initial conditions, \eqref{90} follows, hence the result.
\qed

\subsection{The small ball exit measure}\label{sec34}

This section and the next are devoted to the proof of Proposition \ref{lem3}.
By assumption, the parameters $\la^r_i$ and $\mu^r_i$ scale like $r^2$,
as expressed in equation \eqref{204}. The special case where, for all $r$,
$\la^r_i=\la_ir^2$ and $\mu^r_i=\mu_ir^2$ is referred to as {\it the homogeneous case}.
Our strategy is to first prove the lemma in the homogeneous case,
where the processes $\hat X^r$ can all be
expressed as scaled versions of a {\it single process}.
This is the content of this section.
In \S\ref{sec35}, the general case is considered, and by appealing
to a change of measure argument, Proposition \ref{lem3} is proved.

\begin{lemma}
  \label{lem5}
The statement of Proposition \ref{lem3} holds in the homogeneous case,
with $u(r)=r^{-\kap_0}$, for any $\kap_0\in(0,1/2)$.
\end{lemma}

\noi{\bf Proof.}
Let $2N$ mutually independent Poisson processes
$A_i, S_i$, be given, with intensities $\la_i$ and $\mu_i$, respectively.
Since by assumption $\mu^r_i=\mu_ir^2$ and $\la^r_i=\la_ir^2$,
the tuple $(A^r_i,S^r_i)$ is equal in law, for each $r$, to
$(A_i(r^2\cdot),S_i(r^2\cdot))$, and without loss of generality we may, and will, assume
that $A^r_i(t)=A_i(r^2t)$ and $S^r_i(t)=S_i(r^2t)$ for all $r$, $i$ and $t$.
Let now $Q$, $X$, $T$ and $D$ be defined as the processes $Q^1$, $X^1$, $T^1$ and
respectively $X^1$ (that is, $Q=Q^r$ where one sets $r=1$).
Then in particular, equations \eqref{asaf1}, \eqref{37} and \eqref{38}
are satisfied by $Q$, $X$, $T$ and $D$, and $X$ is a Markov process on
$\calS^1_{\rm u}$ (see \eqref{28-}).

Since for each $r$ we have the aforementioned relation
between $(A^r_i,S^r_i)$ and $(A_i,S_i)$, one can also express
$(Q^r,X^r,T^r,D^r)$ as certain path transformations of $(Q,X,T,D)$,
for each $r$. The most significant aspect of this in the proof is that
the rescaled processes $\hat X^r$ and $\hat R^r$ can be written as
rescaled versions of a single process. Denote $R=\one\cdot X$.
Then by \eqref{asaf1}, $X_i=(\mu_i)^{-1}Q_i$, whereas
$X^r_i=(r^2\mu_i)^{-1}Q^r_i$. Hence $X^r=r^{-2}X(r^2\cdot)$, and thus
by \eqref{14},
\begin{align}\label{101}
\hat X^r(t)=r^{-1}X(r^2t),
\qquad
\hat R^r(t)=r^{-1}R(r^2t).
\end{align}
The state space $\calS^r$ for the Markov process $\hat X^r$
is given in the case under consideration as
$\calS^r=r^{-1}\calS^1_{\rm u}=(r\mu_1)^{-1}\Z_+\times\cdots\times(r\mu_N)^{-1}\Z_+$.

For $\eps>0$, set
\begin{equation}\label{17}
\tau(\eps)=\inf\{t\ge0:R(t)\ge \eps\}.
\end{equation}
Clearly we have the identity $\tau^r(\eps)=r^{-2}\tau(\eps r)$ (see \eqref{35}).
Let $\kap_0$ be as in the statement of the lemma, that is, $\kap_0\in(0,\frac{1}{2})$,
and set $\kap=1-\kap_0\in(\frac{1}{2},1)$. Denote
\[
B^r_i=B(re_i,r^{\kap})
\]
and $q^r=(q^r_i)_i$, where
\[
q^r_i=\PP^1_0(X(\tau(r))\in B^r_i),
\]
where $\PP^1_x$ (with the corresponding expectation $\E^1_x$) stands for the law of $X=X^1$ with $X(0)=x$.
Then
\begin{align}\label{100}
\PP^r_0(\hat X^r(\tau^r(\eps))\in B(\eps e_i,r^{-\kap_0}))
&=\PP^1_0(X(\tau(\eps r))\in B(\eps r e_i, r^{1-\kap_0}))
\\
&\ge \PP^1_0(X(\tau(\eps r))\in B(\eps r e_i, (\eps r)^{1-\kap_0}))
\notag
\\
&\ge\PP^1_0(X(\tau(\eps r))\in B^{\eps r}_i)
\notag
\\
\notag
&=q^{\eps r}_i.
\end{align}
Using the fact that the balls $B(\eps e_i,r^{-\kap_0})$
are disjoint for each $\eps$ and sufficiently large $r$, \eqref{200}
will follow once we show that
\begin{align}\label{15}
\text{there exists $q\in\calM_1$ such that $\lim_{r\to\iy}q^r=q$}.
\end{align}
Note that Proposition \ref{lem3} asserts, moreover, that $q$ does not depend on the
choice of $\kap_0$. To address this point,
consider $0<\kap_0<\kap_0'<\frac{1}{2}$ for which $q$ and, respectively $q'$,
satisfy \eqref{15}. Then the fact that the LHS of \eqref{100}
is monotone decreasing in $\kap_0$ gives $q_i\ge q'_i$ for all $i$,
and since $q$ and $q'$ are members of $\calM_1$, this shows that $q=q'$.
Hence the proof will be complete once we show \eqref{15} for fixed $\kap_0$.

To this end, note that it suffices to show
that there exist $\delta\in(0,1)$ and $K>0$ such that for every
$k\in\N$, $k\ge K$, one has
\begin{equation}
  \label{16}
  |q^r_i-q^m_i|\le\delta^k,
  \text{ for all $i\in[N]$ and $r\in[2^k,2^{k+1}]$, where $m=2^{k+2}$}
\end{equation}
and
\begin{align}\label{25}
\lim_{r\to\iy}\sum_iq^r_i=1.
\end{align}
Indeed, if $r$ and $m$ are both within $[2^k,2^{k+1}]$
then \eqref{16} gives $|q^r_i-q^m_i|\le2\delta^k$. As a result, for arbitrary
$r<m$, denoting $a(\ell)=\lfloor\log_2(\ell)\rfloor$,
\[
|q^r_i-q^m_i|\le\sum_{j=a(r)}^{a(m)}2\delta^j\le 2(1-\delta)^{-1}\delta^{a(r)}.
\]
This shows that, for fixed $i$, any sequence $\{q^r_i\}_r$ is a Cauchy sequence as $r\to\iy$.
Along with \eqref{25}, we obtain that \eqref{15} holds.

In what follows, we prove \eqref{16} and \eqref{25}.
We let $r\in[2^k,2^{k+1}]$ and $m=2^{k+2}$, where $k$ is arbitrary, but fixed.
Without loss of generality, the proof of \eqref{16} considers only $i=1$.

\begin{figure}
\begin{center}
\includegraphics[width=16em,height=13em]{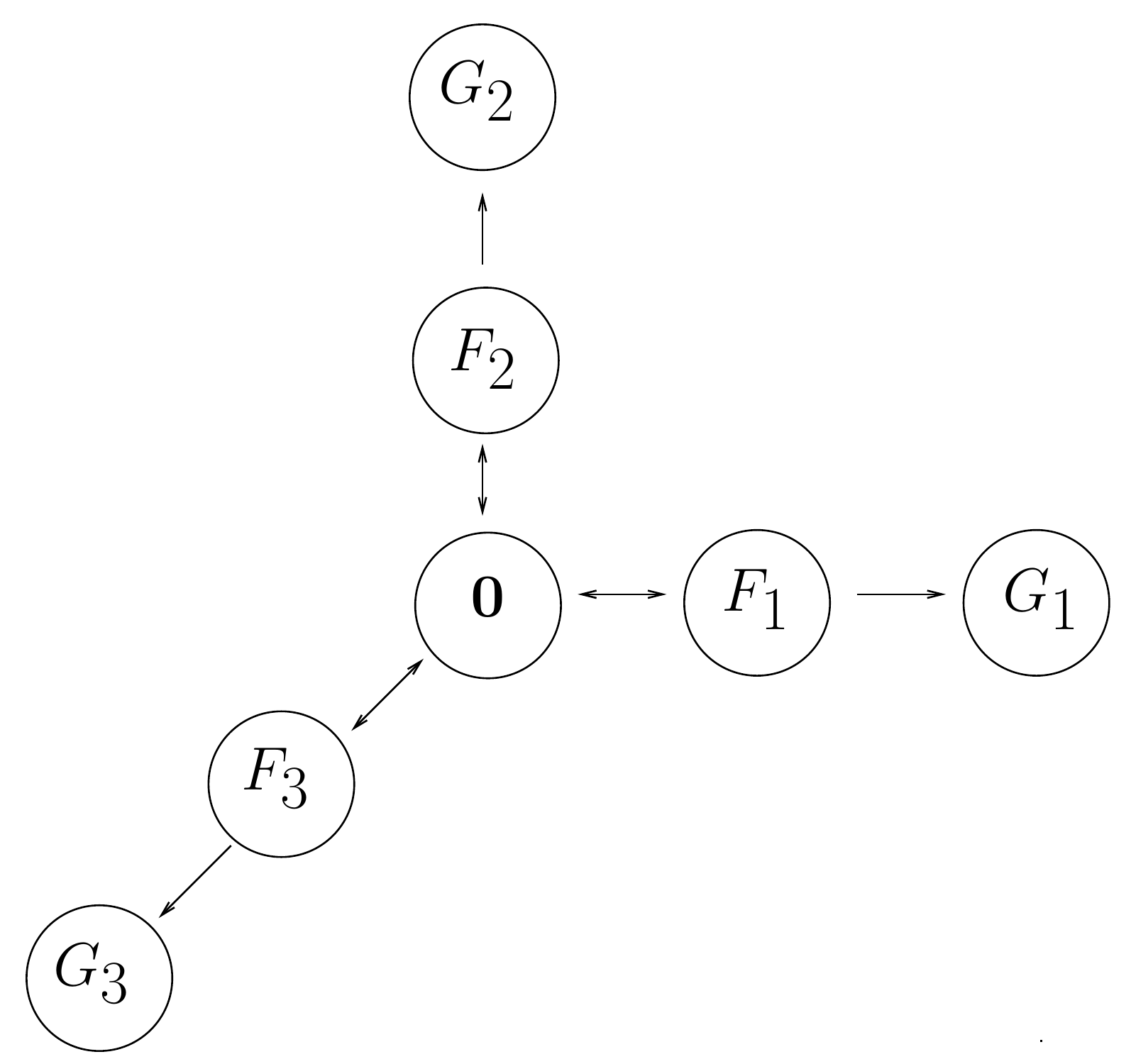}
\caption{\sl\footnotesize
In this toy model, denoting transition probabilities by $p(\cdot,\cdot)$, the quantity $\max_{i,j}|p(F_i,0)-p(F_j,0)|$
controls
$\max_i|p(0,F_i)-P_0(\text{the process is absorbed at } G_i)|$.
}
\end{center}
\end{figure}

To help explain the main idea and motivate a couple of technical tools,
we first consider a highly simplified model, illustrated in Figure 1.
Consider a discrete time Markov process on a finite set $S$ that is star shaped. That is,
$S$ consists of $2N+1$ states, denoted by $0$, $F_i$, $G_i$, $i\in[N]$.
For each $i$, $F_i$ communicates only with $0$ and $G_i$.
The state $0$ communicates only with the states $F_i$,
while $G_i$ are absorbing. Denoting transition probabilities
by $p(s,s')$, we have $p(0,F_i)>0$, $p(F_i,0)>0$,
$p(F_i,G_i)>0$, and $p(G_i,G_i)=1$, while all other transition probabilities are zero.
For $s\in S$, let $\bar p(s,G_i)$ denote the probability to get absorbed at $G_i$
starting from $s$. Then
\begin{align*}
  \bar p(0,G_1)&=\sum_{i\ge1} p(0,F_i)\bar p(F_i,G_1)\\
  &=p(0,F_1)p(F_1,G_1)+\sum_{i\ge1} p(0,F_i)p(F_i,0)\bar p(0,G_1).
\end{align*}
From this one obtains
\[
\frac{\bar p(0,G_1)}{p(0,F_1)}=\frac{1-p(F_1,0)}{1-\sum_{i\ge1}p(0,F_i)p(F_i,0)}.
\]
If the transition probabilities starting at $F_i$ depend weakly on $i$,
in the sense that for some $\delta>0$ one has
$|p(F_i,0)-p(F_1,0)|<\del$ for all $i$, and if in addition
$1-p(F_1,0)>c>0$ for some constant $c$, it follows that
\[
|\bar p(0,G_1)-p(0,F_1)|\le c'\del,
\]
where $c'$ depends on $c$ but not on $\del$.
The relevance to our problem is as follows.
Roughly speaking, the states $F_i$ and $G_i$ represent the collections
of states within $B^r_i$ and $B^m_i$, respectively, and $0$ represents the origin.
The calculation above suggests that if
the probability of reaching $0$ before reaching $B^m_i$
starting anywhere in $B^r_i$ depends weakly on $i$
then the difference $p^r_i-p^m_i$ is small.

We now consider the process $X$ and the stopping times $\tau(\eps)$,
and in addition let
\[
\zeta=\inf\{t\ge0:X(t)=0\}.
\]
We aim at showing there exist $r_0,c>0$ such that for every $r>r_0$, one has
\begin{equation}
  \label{20}\text{for every } i\in[N] \text{ and } x\in B^r_i,
  \ \Big|\PP^1_x(\zeta<\tau(m))-\frac{m-r}{m}\Big|\le r^{-c},
\end{equation}
\begin{equation}
  \label{24}
  \PP^1_0(r^{-2}\tau(r)>t)\le e^{-ct},
\end{equation}
\begin{equation}
  \label{18}
  \PP^1_0(X(\tau(r))\notin\cup_i B^r_i)\le r^{-c},
\end{equation}
\begin{equation}
  \label{19}
  \text{for every } i\in[N] \text{ and } x\in B^r_i,
  \ \PP^1_x(\tau(m)<\zeta,\, X(\tau(m))\notin B^m_i)\le r^{-c}.
\end{equation}
For estimate \eqref{20}, note that it follows from the identity $R=\one\cdot X$,
the expression \eqref{10} for the generator of $X$ (with $\la_i$ and $\mu_i$
substituted for $\la^r_i$ and $\mu^r_i$), and condition \eqref{11},
that the stopped process $R(\cdot\wedge\zeta\wedge\tau(m))$ is a martingale.
By this martingale property and the fact that $X$ lives on the grid $\calS^1_{\rm u}$,
there is a constant $c$ such that
\begin{align}\notag
\frac{m-r-r^\kap-c}{m+c}
\le\PP^1_x(\zeta<\tau(m))\le\frac{m-r+r^\kap+c}{m}.
\end{align}
Estimate \eqref{20} follows, using the fact that $m\ge 2r$.

For inequality \eqref{24}, the relation $r^{-2}\tau(r)=\tau^r(1)$ gives
$\PP^1_0(r^{-2}\tau(r)>1)=\PP^r_0(\|\hat R^r\|_1<1)$.
By Lemma \ref{lem1}, $\hat R^r$ converges in law to an RBM $\rho$, and so
$\PP^r_0(\|\hat R^r\|_1<1)\to \PP_0(\|\rho\|_1<1)<1$.
Thus there exists $\gamma\in(0,1)$ such that for all $r$ sufficiently large,
$\PP^1_0(r^{-2}\tau(r)>1)\le \gamma$. For other initial conditions $x$
the probability of this event under $\PP^1_x$ is even smaller, and
therefore is still bounded by $\gamma$. Markovity thus gives \eqref{24}.

For estimate \eqref{18}, fix $c_0>0$ to be a constant $c$ that satisfies \eqref{24}.
Recall the definition of $F$ from \eqref{eq100}.
Since $X$ takes values in $\calS^1_{\rm u}$, $\one\cdot X(\tau(r))$
must take a value within $[r,r+c_1]$, for some constant $c_1>0$.
We claim that if $x\in\R_+^N$, $\one\cdot x\in[r,r+c_1]$
and $x\notin\cup_iB^r_i$ then $F(x)\ge r^{\kap}/\sqrt{2}$.
To this end, assume, without loss of generality, that $x_1=\max_ix_i$.
Since $x\notin B^r_1$, we have
$(x_1-r)^2+\sum_{i=2}^Nx^2_i\ge r^{2\kap}$. Therefore,
if $r-r^\kap/\sqrt{2}\le x_1\le r+c_1$, we have
that $F(x)^2=(\sum_{i=2}^Nx_i)^2\ge r^{2\kap}/2$
and hence $F(x)\ge r^{\kap}/\sqrt{2}$.
If on the other hand $x_1<r-r^\kap/\sqrt{2}$, then using $\one\cdot x\ge r$,
we obtain again $F(x)=\sum_{i=2}^Nx_i>r^\kap/\sqrt{2}$.

As a result of the above claim, for all large $r$,
\begin{align}\label{202}
  \PP^1_0(X(\tau(r))\notin\cup_i B^r_i)
  &\le
  \PP^1_0(\tau(r)>c_0r^2\log r)+\PP^1_0(\|F(X(\cdot))\|_{c_0r^2\log r}\ge r^\kap/\sqrt{2}).
\end{align}
By \eqref{24}, the first term above is bounded by $r^{-c_0}$.
Since $F(X(0))=0$, we have by Lemma \ref{lem2}, relation \eqref{101},
and the relation $\kap=1-\kap_0$, that the second term is bounded by $r^{-c}$.

For estimate \eqref{19}, define
\begin{align}\notag
\nu(m)=\inf\{t\ge0 : R(t)\le m^\kap\}.
\end{align}
Fix $i$ and consider $x\in B^r_i$. Denote $B_{0,m}=B(0,m^\kap)\cap\calS^1_{\rm u}$.
Then a use of Strong Markovity gives
\begin{equation}\label{102}
\PP^1_x(\tau(m)<\zeta, X(\tau(m))\notin B^m_i)
\le\PP^1_x(\tau(m)<\nu(m), X(\tau(m))\notin B^m_i)
+\max_{z\in B_{0,m}}\PP^1_z(\tau(m)<\zeta).
\end{equation}
To bound the first term
consider the event $\tau(m)<\nu(m)$.
If $X(\tau(m))\in\cup_{j\ne i}B^m_j$
then $\|F(X(\cdot))\|_{\tau(m)}>m^\kap$ holds,
whereas if $X(\tau(m))\notin\cup_{j}B^m_j$, then by the argument
provided in the previous paragraph
we have $\|F(X(\cdot))\|_{\tau(m)}\ge m^\kap/\sqrt{2}$.
This implies that the first term in \eqref{102} is bounded by
\[
\PP^1_x(\tau(m)\wedge\nu(m)>c_0m^2\log m)
+\PP^1_x(\|F(X(\cdot))\|_{c_0m^2\log m}>m^\kap/\sqrt{2}).
\]
This expression can be handled as the RHS of \eqref{202}.
Since $x\in B^r_i$, $\|F(X(0))\|\le r^\kap
\le \frac{1}{2^\kap}m^\kap$. Thus Lemma \ref{lem2}
is applicable with $\gamma_1=2^{-\kap}<\gamma_2=1/\sqrt{2}$.
This gives the bound $r^{-c}$ on the first term on the RHS of \eqref{102}.

To bound the second term in \eqref{102} we again use the martingale property of
$R(\cdot\wedge\zeta\wedge\tau(m))$. It gives
\begin{align}\notag
\PP^1_z(\tau(m)<\zeta)\le\frac{m^\kap}{m-2m^\kap}.
\end{align}
Since $2r\le m\le 4r$, we obtain that for sufficiently large $r$,
the last term in \eqref{102} is bounded above by $r^{-c}$.
This completes the proof of \eqref{20}--\eqref{19}.

We now deduce \eqref{16} and \eqref{25} from \eqref{20}--\eqref{19}.
Identity \eqref{25} follows immediately from \eqref{18}.
For $x\in\calS^r_{\rm u}$ (see \eqref{28-}) denote
\begin{align}\notag
q(x,r,m)=\PP^1_0(X(\tau(r))=x)\PP^1_x(X(\tau(m))\in B^m_1)
\end{align}
and $B_{r,i}=B^r_i\cap\calS^r_{\rm u}$. Then
\begin{align}
q^m_1&=
\label{21}
\sum_{x\in \calS^r_{\rm u}:x\notin\cup_{i\ge1}B^r_i}q(x,r,m)
+\sum_{i>1}\sum_{x\in B_{r,i}}q(x,r,m)
+\sum_{x\in B_{r,1}}q(x,r,m)\\
\notag
&=:\beta^{r,m,1}+\beta^{r,m,2}+\beta^{r,m,3}.
\end{align}
It follows from \eqref{18} that $\beta^{r,m,1}\le r^{-c}$.
Next, consider the term $\beta^{r,m,2}$. Let $i>1$ and $x\in B_{r,i}$.
Then
\[
\PP^1_x(X(\tau(m))\in B^m_1)
=\PP^1_x(\zeta<\tau(m),X(\zeta(m))\in B^m_1)+\PP^1_x(\zeta(m)<\zeta,X(\tau(m))\in B^m_1).
\]
The first term above is equal to $\PP^1_x(\zeta<\tau(m))q^m_1$.
By \eqref{19}, the second term bounded by $r^{-c}$.
Combining this with \eqref{20},
\[
\beta^{r,m,2}=\sum_{i>1}q^r_i\frac{m-r}{m}q^m_1+\eps(r,m),
\]
where here and in the remainder of this proof, $\eps(r,m)$
denotes a generic function of $(r,m)$ which staisfies
$|\eps(r,m)|\le r^{-c}$ for all large $r$.

As for $\beta^{r,m,3}$, consider $x\in B^r_1$. We have
\begin{align*}
&\PP^1_x(X(\tau(m))\in B^m_1)
\\
&\quad=\PP^1_x(\zeta<\tau(m),X(\tau(m))\in B^m_1)
+\PP^1_x(\tau(m)<\zeta,X(\tau(m))\in B^m_1)\\
&\quad=\PP^1_x(\zeta<\tau(m))q^m_1+\PP^1_x(\tau(m)<\zeta)
-\PP^1_x(\tau(m)<\zeta,X(\tau(m))\notin B^m_1).
\end{align*}
Using \eqref{19}, the last term above is bounded, in absolute value, by $r^{-c}$.
Combined with \eqref{20}, this gives
\[
\beta^{r,m,3}=q^r_1\Big(\frac{m-r}{m}q^m_1+\frac{r}{m}\Big)+\eps(r,m).
\]
Combining the three estimates,
\begin{align*}
q^m_1&=\sum_{i=1}^rq^r_i\frac{m-r}{m}q^m_1+q^r_1\frac{r}{m}+\eps(r,m)\\
&=(1-\eps(r,m))\frac{m-r}{m}q^m_1+q^r_1\frac{r}{m}+\eps(r,m),
\end{align*}
where \eqref{18} is used. Hence, using $m/r\le 4$,
\[
|q^m_1-q^r_1|\le\frac{m}{r}|\eps(r,m)|\le cr^{-c}.
\]
This gives \eqref{16} and completes the proof of the lemma.
\qed

\subsection{Relaxation of the homogeneity assumption}\label{sec35}

In this section we prove Proposition \ref{lem3} based on Lemma \ref{lem5},
by means of a change of measure.
Thus the general setting, where $\la^r_i$ and $\mu^r_i$
satisfy the hypotheses of Theorem \ref{th1}, is in force.
Since the statement of Proposition \ref{lem3} refers to $\PP^r_0$,
we may and will assume in this section that the initial condition is
$Q^r(0)=X^r(0)=0$ identically. Thus the only stochastic primitives
in the model are the processes $(A^r,S^r)$. In particular,
as follows from equations \eqref{asaf1}, \eqref{37}, \eqref{38-}, \eqref{38} and \eqref{14},
the processes $X^r(t)$, $t\in[0,T]$ and $\hat X^r(t)$, $t\in[0,T]$
are determined by $(A^r(t),S^r(t))$, for $t\in[0,T]$.
In addition to the measure $\bP$, we introduce below
a reference probability measure
$\bQ$ on $(\Om,\calF)$ under which, for all $r$, the Poisson processes $A^r_i$ and $S^r_i$
have intensities $\la^{0,r}_i$ and $\mu^{0,r}_i$, respectively,
where we denote $\la^{0,r}_i=\la_ir^2$ and $\mu^{0,r}_i=\mu_ir^2$.
Denote by $\E_{\bQ}$ the corresponding expectation.
The laws of the driving Poisson processes as well as
that of the queue length
process $Q^r$ under $\bP$ can then be obtained from those under $\bQ$
by a change of measure (as shown below). However, this does not apply to
the nominal workload process $X^r$, for which the parameters $\la^r_i$ and $\mu^r_i$
determine not only the jump intensities
but also the scaling factors in the definition
\eqref{asaf1} of $X^r$ in terms of $Q^r$. This is reflected also
in the formula for the generator
$\calL^r_{\rm u}$ (see \eqref{10}) where these parameters enter in both
the jump rates and the jump sizes.
An intermediate transformation is required.

To this end we define analogously to \eqref{asaf1} and \eqref{14},
a process $X^{0,r}_i$ and its scaled version by
\begin{align}\notag
X^{0,r}_i=(\mu^{0,r}_i)^{-1}Q^r_i,
\qquad
\hat X^{0,r}_i=rX^{0,r}_i.
\end{align}
Similarly, we let $\hat R^{0,r}=\one\cdot\hat X^{0,r}$
and $\tau^{0,r}(\eps)=\inf\{t\ge0:\hat R^{0,r}(t)\ge\eps\}$.

The starting point of this section is to notice
that Lemma \ref{lem5}, proved in the previous
section, implies that there exists $q\in\calM_1$ such that for $\kap_0\in(0,1/2)$,
\begin{equation}
  \label{103}
  \lim_{\eps\down0}\limsup_{r\to\iy}
  |\bQ(\hat X^{0,r}(\tau^{0,r}(\eps))\in B(\eps e_i,r^{-\kap_0}))-q_i|=0,
  \qquad i\in[N].
\end{equation}
The proof proceeds in two steps. First, it is shown that
a version of \eqref{103}, that refers to $\hat X^r(\tau^r(\eps))$ in place of
$\hat X^{0,r}(\tau^{0,r}(\eps))$, is valid,
and then that the same statement remains true under $\bP$ (equivalently,
under $\PP_0^r$).

\skp

\noi{\bf Proof of Proposition \ref{lem3}.}
We first prove that, for $q$ as in \eqref{103}, there exists $u\in\calU_0$ such that
\begin{equation}
  \label{104}
  \lim_{\eps\down0}\limsup_{r\to\iy}
  |\bQ(\hat X^r(\tau^r(\eps))\in B(\eps e_i,u(r)))-q_i|=0,
  \qquad i\in[N].
\end{equation}
Based on \eqref{103}, the statement \eqref{104}
is almost an immediate consequence of convergence of $\hat R^{0,r}$ to an RBM
under $\bQ$ and the closeness of $\hat X^{0,r}$ and $\hat X^r$.
Indeed, the relation between $\hat X^{0,r}$ and $\hat X^r$ is
$\hat X^{0,r}_i=\beta^r_i\hat X^r_i$ where $\beta^r_i=\mu^r_i/\mu^{0,r}_i$.
We have $\max_i|\beta^r_i-1|<cr^{-1}$ by \eqref{204}.
Thus for $\eps<1$, $\|\hat X^{0,r}(t)-\hat X^r(t)\|<cr^{-1}$ for all $t\le\tau^r(\eps)
\w\tau^{0,r}(\eps)$. Hence \eqref{104} will follow from \eqref{103}
if we show that, as $r\to\iy$,
$\sup_{\eps\in(0,1)}\|\hat X^{0,r}(\tau^r(\eps))-\hat X^{0,r}(\tau^{0,r}(\eps))\|\to0$
in probability.
Since $\hat X^{0,r}$ are $\calC$-tight by Lemma \ref{lem-b} and $\tau^{0,r}(\eps)$
are dominated by $\tau^{0,r}(1)$, that form a tight sequence of RVs,
it suffices to prove that
\begin{equation}\label{105}
\sup_{\eps\in(0,1)}|\tau^r(\eps)-\tau^{0,r}(\eps)|\to0 \text{ in probability}.
\end{equation}
The convergence of $\hat R^{0,r}$ to RBM implies that for any $M>0$ and $\del>0$,
\[
\lim_{\kap\down0}\limsup_{r\to\iy}
\bQ\Big(\inf_{t\in[0,M]}\sup_{u\in(0,\del)}
|\hat R^{0,r}(t+u)-\hat R^{0,r}(t)|<\kap\Big)=0.
\]
It follows that for any $M>0$ and $\del>0$,
\[
\limsup_{r\to\iy}\bQ\Big(\sup_{\eps\in(0,1)}|\tau^r(\eps)\w M-\tau^{0,r}(\eps)\w M|>\del\Big)=0.
\]
Using again the tightness of the RVs $\tau^{0,r}(1)$, \eqref{105} follows,
hence also \eqref{104}.

The second and final step is to prove that in \eqref{104},
$\bQ$ may be replaced by $\bP$.
Denote the events of interest by $K^r_{\eps,i}=
\{\hat X^r(\tau^r(\eps))\in B(\eps e_i,u(r))\}$.
Since $\sum_iq_i=1$, using the fact that for any $\eps$ and all
sufficiently large $r$, $\{K^r_{\eps,i}\}_{i}$ are disjoint,
it suffices to prove for each $i$ the lower bound
\begin{equation}\label{108}
\liminf_{\eps\down0}\;\liminf_{r\to\iy}\bP(K^r_{\eps,i})\ge q_i.
\end{equation}
Given any $\del>0$, we clearly have
$\lim_{\eps\down0}\limsup_{r\to\iy}\bQ(\tau^r(\eps)>\del)=0$.
Hence by \eqref{104}, denoting
$K^r_{\eps,\del,i}=K^r_{\eps,i}\cap\{\tau^r(\eps)\le\del\}$,
we have
\begin{equation}
  \label{106}
  \lim_{\eps\down0}\limsup_{r\to\iy}
  |\bQ(K^r_{\eps,\del,i})-q_i|=0,
  \qquad i\in[N].
\end{equation}

A change of measure is formulated in terms of the exponential martingale
\begin{align}\notag
\psi^r_t=\exp\sum_i\Big[A^r_i(t)\log\Big(\frac{\la^r_i}{\la^{0,r}_i}\Big)
- (\la^r_i-\la^{0,r}_i )t
+S^r_i(t)\log\Big(\frac{\mu^r_i}{\mu^{0,r}_i}\Big)
- (\mu^r_i-\mu^{0,r}_i )t\Big].
\end{align}
Let $\calA^r_t=(A^r_i(u),S^r_i(u))_{i\in[N],u\in[0,t]}$.
Let also $\calG^r_t=\sig\{\calA^r_t\}$. For each $r$ and $t$,
let a probability measure $\bP^{r,t}$ on $(\Om,\calG^r_t)$ be defined by
$\bP^{r,t}(G)=\E_{\bQ}[\psi^r_t\one_G]$ for $G\in\calG^r_t$.
Then, for each $r$ and $t$, the law of $\calA^r_t$ under $\bP^{r,t}$
is the same as that under $\bP$.
Moreover, note that for each $i$,
the event $K^r_{\eps,\del,i}$ is measurable on $\calG^r_\del$.
Hence to establish \eqref{108}, it suffices to prove that for each $i$,
\begin{equation}\label{109}
\tilde q_i:=\liminf_{\del\down0}\;
\liminf_{\eps\down0}\;\liminf_{r\to\iy}\E_{\bQ}[\psi^r_\del\one_{K^r_{\del,\eps,i}}]
\ge q_i.
\end{equation}

For $\eta>0$, denote
$G^r_{\del,\eta}=\{\psi^r_\del>1-\eta\}$.
Suppose we show that for any $\eta>0$,
\begin{equation}\label{107}
\liminf_{\del\down0}\;\liminf_{r\to\iy}\bQ(G^r_{\del,\eta})=1.
\end{equation}
Then we may argue as follows,
\[
\E_{\bQ}[\psi^r_\del\one_{K^r_{\eps,\del,i}}]
\ge\E_{\bQ}[\psi^r_\del\one_{K^r_{\eps,\del,i}\cap G^r_{\del,\eta}}]
\ge(1-\eta)\bQ(K^r_{\eps,\del,i})-\bQ((G^r_{\del,\eta})^c).
\]
Taking $r\to\iy$ then $\eps\down0$, then using \eqref{106}, and finally
taking $\del\down0$, gives
\[
\tilde q_i\ge(1-\eta)q_i-\limsup_{\del\down0}\;\limsup_{r\to\iy}\bQ((G^r_{\del,\eta})^c)
=(1-\eta)q_i.
\]
Since $\eta>0$ is arbitrary, this gives \eqref{109}
and consequently \eqref{108}.

Thus the proof will be complete once \eqref{107} is shown.
To this end, let
\[
\tilde A^r_i(t)=r^{-1}(A^r_i(t)-\la^{0,r}_it),
\qquad
\tilde S^r_i(t)=r^{-1}(S^r_i(t)-\mu^{0,r}_it).
\]
These processes, defined analogously to $\hat A^r$ and $\hat S^r$,
converge under $\bQ$ to BMs. Denote $\hat\la^r_i=r^{-1}(\la^r_i-\la^{0,r}_i)$
and $\hat\mu^r_i=r^{-1}(\mu^r_i-\mu^{0,r}_i)$ and recall
by \eqref{204} that these sequences converge.
Write $\psi^r_t$ in terms of $\tilde A^r$ and $\tilde S^r$ as
\[
\psi^r_t=\exp\sum_i\Big[\tilde A^r_i(t)rU^r_i+(\la^{0,r}_iU^r_i-r\hat\la^r_i)t
+\tilde S^r_i(t)rV^r_i+(\mu^{0,r}_iV^r_i-r\hat\mu^r_i)t\Big],
\]
where $U^r_i=\log(1+\frac{\hat\la^r_i}{r\la_i})$,
$V^r_i=\log(1+\frac{\hat\mu^r_i}{r\mu_i})$.
Denoting $L^r_t=\max_i|\tilde A^r_i(t)|\vee|\tilde S^r_i(t)|$
and using $|\log(1+x)-x|\le cx^2$ for all $|x|<1/2$,
we have for all large $r$,
\[
\log\psi^r_t\ge -cL^r_t-ct.
\]
The aforementioned convergence to BM clearly implies that, for any $\eta>0$,
\[
\liminf_{\del\down0}\liminf_{r\to\iy}\bQ(L^r_\del<\eta)=1.
\]
We thus obtain \eqref{107} and complete the proof.
\qed

\section{Concluding remarks}\label{sec4}

\begin{enumerate}

\item
It is desirable to extend the main result of this paper
beyond the Markovian setting, to general service time
distributions and renewal arrival distributions, under second
moment conditions.
Whereas the behavior of the modulus according to an RBM
certainly holds in vast generality,
and the attraction to the collection of axes $\calS_0$ can likely be
extended, the existence of limiting
entrance laws appears to require different machinery.
Indeed, the proof presented here makes crucial
use of the strong Markovity of the prelimit processes.

\item
The proof presented in this paper sheds no light on
the angular distribution $q$ (except that
it does not depend on the second order parameters $\hat\la_i$, $\hat\mu_i$).
A characterization of $q$ that would be useful and lead to
further information about it is desirable.

\item
Figure 2 depicts results of Monte Carlo simulations
for an SSQ model with $N=2$ at criticality, aimed at estimating $q$.
It shows the behavior of $q_1$ as several parameters vary.
They all suggest monotone dependence, that one would wish
to substantiate mathematically.
\begin{enumerate}
\item
The graph shown at Figure 2(c) is, in particular, relevant
to the heuristic mentioned in the introduction,
according to which more variable traffic attains lower priority.
In this example, the traffic intensities $\la_i/\mu_i$ are kept fixed.
As $\la_1$ increases, the inter-arrival variance increases,
which, according to the graph, increases $q_1$,
indicating lower priority for this class.
\item Figure 2(d) shows the dependence on the tie breaking
parameter, $p_1$. It exhibits that
the tie breaking rule affects the limiting angular distribution.
However, we have not aimed at providing a proof of this claim.
\end{enumerate}

\end{enumerate}

\begin{figure}
\footnotesize
\begin{center}
\includegraphics[width=20em,height=11em]{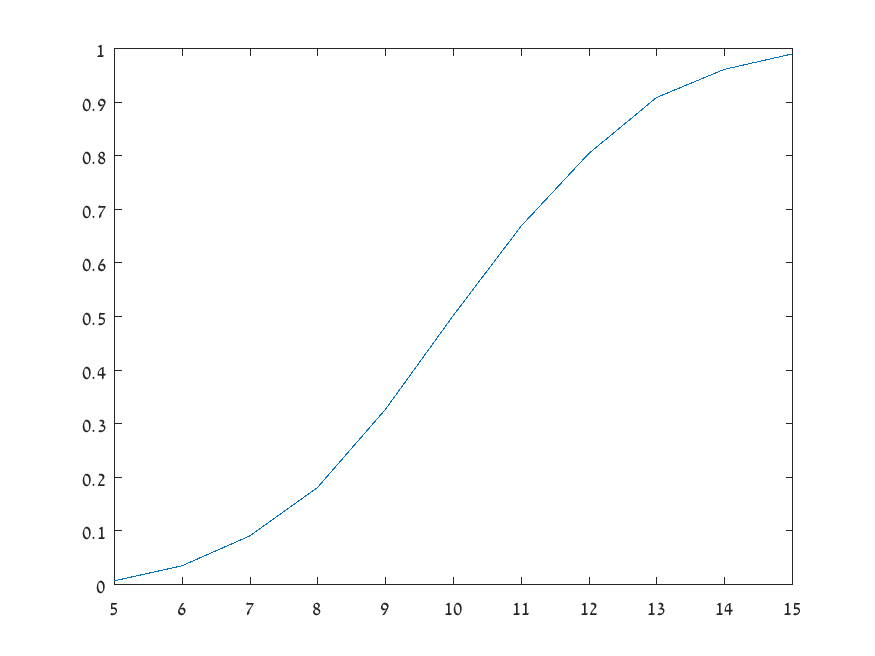}
\hspace{3em}
\includegraphics[width=20em,height=11em]{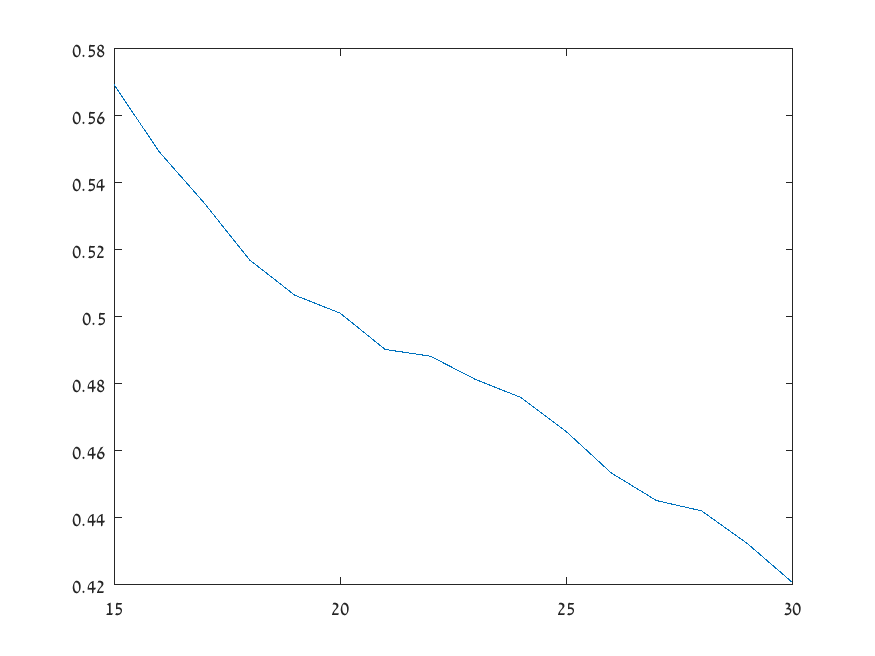}
\\ \vspace{-0.9em}
$\la_1$
\hspace{22em}
$\mu_1$
\\
(a) \hspace{21.6em} (b)
\\ \ \\
\includegraphics[width=20em,height=11em]{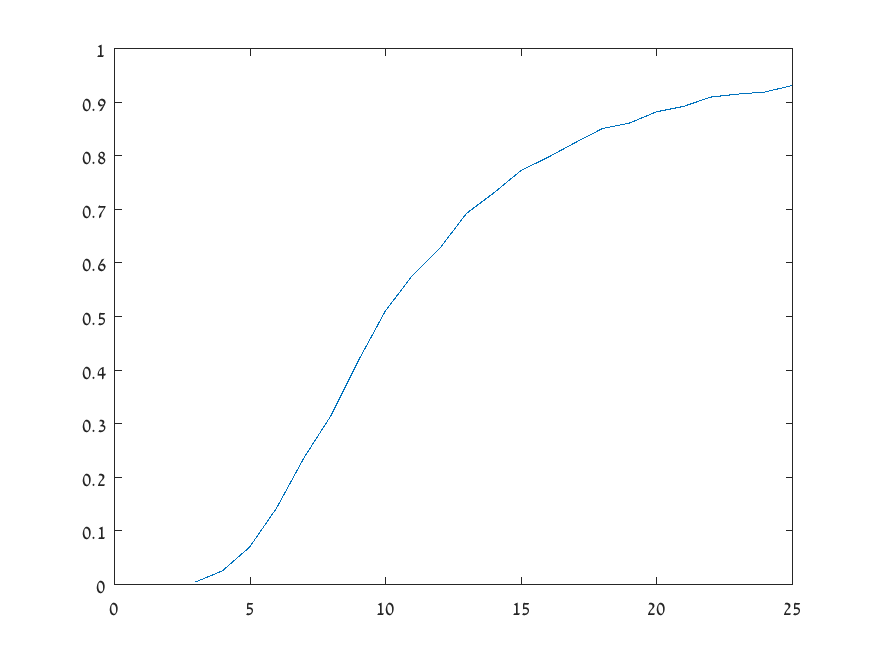}
\hspace{3em}
\includegraphics[width=20em,height=11em]{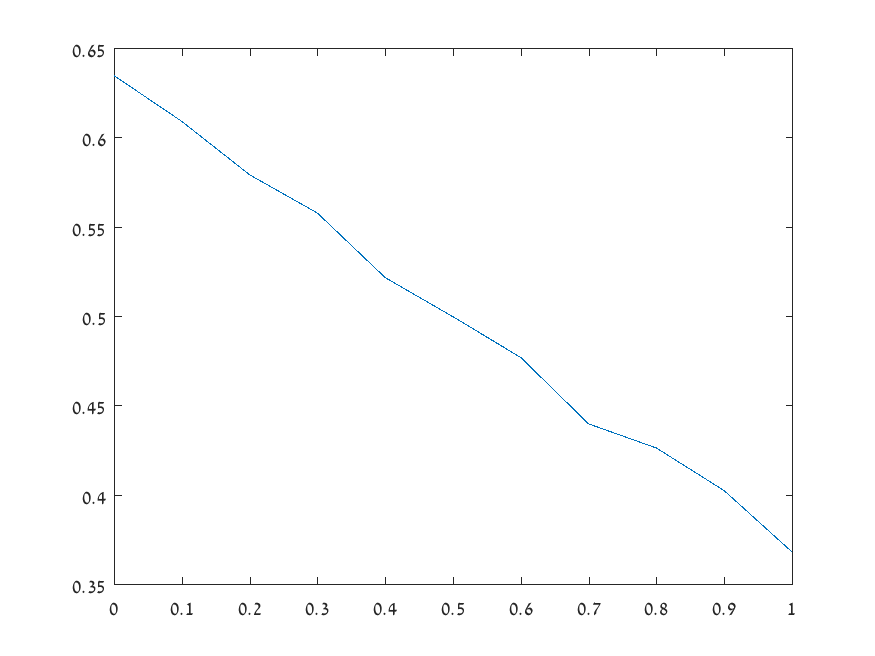}
\\ \vspace{-0.9em}
$\la_1$
\hspace{22em}
$p_1$
\\
(c) \hspace{21.6em} (d)
\caption{\sl\footnotesize
Simulation of $q_1$ as a function of various parameters, for $N=2$.\newline
(a) $q_1$ as a function of $\la_1$, fixed $\mu$'s:
$\mu_1=\mu_2=20$, $\la_2=\mu_2-\la_1$.\newline
(b) $q_1$ as a function of $\mu_1$, fixed $\la$'s:
$\la_1=\la_2=10$, $\mu_2=1/(1/\la_1-1/\mu_1)$.\newline
(c) $q_1$ as a function of $\la_1$, fixed ratio $\la_1/\mu_1$, $\la_2$ and $\mu_2$:
$\la_2=10$, $\mu_2=20$, $\mu_1=2\la_1$.\newline
(d) $q_1$ as a function of $p_1$, fixed $\la$'s and $\mu$'s:
$\la_1=\la_2=10$, $\mu_1=\mu_2=20$, $p_2 =1-p_1$.
}
\end{center}
\end{figure}

\skp

\noi{\bf Acknowledgement.}
The authors are grateful to Ross Pinsky for useful discussions about
heat equation estimates, to Bert Zwart for bringing reference \cite{lamsim} to their attention,
and to the referee for careful reading and valuable comments.

\footnotesize

\bibliographystyle{abbrv}

\bibliography{refs}

\end{document}